\documentclass[a4paper, 11pt]{article}

\usepackage{natbib}
\usepackage{a4wide}
\usepackage{amsmath}
\usepackage{amsthm}
\usepackage{amsfonts}
\usepackage{amssymb}
\usepackage{latexsym}
\usepackage[dvips]{epsfig}
\usepackage{exscale}
\usepackage{psfrag}
\usepackage[final]{fixme}

\usepackage{placeins}           
\usepackage{afterpage}          

\newcommand{\bino}[2]{\left(\begin{array}{c}#1\\#2 \end{array}\right)}

\newcommand{\dx}{\,\text{d$x$}}

\newcommand{\ds}{\, \text{\rm d$s$}}

\newcommand{\del}{\delta}
\newcommand{\dlt}{\delta_1}

\newcommand{\vf}[1]{{\bf #1}} 

\newcommand{\e}{{\rm e}}

\newcommand{\R}{\mathbb{R}}

\newcommand{\N}{\mathbb{N}}

\newlength{\boxwidth}
\newlength{\proofwidth}
\newtheorem{theorem}{Theorem}[section]
\newtheorem{corollary}[theorem]{Corollary}
\newtheorem{lemma}[theorem]{Lemma}
\newtheorem{proposition}[theorem]{Proposition}

\theoremstyle{remark}
\newtheorem*{remark}{Remark}

\theoremstyle{definition}

\usepackage{natbib}
\bibpunct{(}{)}{;}{a}{,}{,}

\begin{document}

\title{Analysis of Linear Difference Schemes
in the Sparse Grid Combination Technique}
\author{
Christoph Reisinger \\ Mathematical Institute, University of Oxford \\
24--29 St Giles,
Oxford, OX1 3LB,
United Kingdom \\
christoph.reisinger@maths.ox.ac.uk}
\maketitle

\begin{abstract}
Sparse grids \citep{zenger:90,bungartz-griebel:04} are tailored to the approximation of smooth high-dimensional
functions. On a $d$-dimensional tensor product space, the number of grid points is
$N = \mathcal O(h^{-1} |\log h|^{d-1})$, where $h$ is a mesh parameter.
The so-called combination technique, based on hierarchical decomposition and extrapolation,
requires specific multivariate error expansions
of the discretisation error on Cartesian grids to hold.
We derive such error expansions for linear difference schemes
through an error correction technique of semi-discretisations.
We obtain overall error formulae of the type $\varepsilon = \mathcal{O} (h^p |\log
h|^{d-1})$
and analyse the convergence, with its dependence on dimension and smoothness, by 
examples of linear elliptic and parabolic problems, with numerical illustrations in up to eight dimensions.
\end{abstract}

\noindent
{\bf Key words.} Sparse grids, combination technique,
error expansions, finite difference schemes

\noindent
{\bf AMS subject classifications.} 65N06, 65N12, 65N15, 65N40


\pagebreak

\section{Introduction}


Models with high-dimensional state spaces play an important role in various applications. 
Important examples arise in financial engineering, more specifically in derivative pricing, which typically
involve computation of expected pay-off functions with respect to a number of
underlying stochastic factors, for instance share prices or interest rates.
In the financial industry, this expectation is most commonly estimated by stochastic simulation,
especially if the dimension exceeds three as is often the case. 
Recent computational studies suggest that sparse grid cubature 
\citep{gerstner-griebel-holtz:09},
and the numerical solution of the
corresponding Feynman-Kac PDE on sparse grids 
\citep{hilber-matache-schwab:05, reisinger:06,leentvaar-oosterlee:08}
are promising alternatives.

To give a further example of high-dimensional spaces,
in quantum mechanics, high-dimensional eigenvalue problems for the
Schr\"odinger equation govern the density functions associated with molecular systems. Again, recent analytical and numerical results
\citep{ garcke-griebel:00,griebel-hamaekers:07,yserentant:10,bachmayr:10,zeiser:10} demonstrate the potential of sparse grids to make decisive progress
in this area.

However, several of these works also give an account of the difficulties encountered, notably for problems with reduced regularity, but even for infinitely smooth solutions in higher dimensions.
Although the asymptotic order of complexity of the sparse grid solution is only weakly dependent on the dimension, the ``constants'' in the error formulae reflect the dependence on both the dimension and on a measure of the variation of the solution. It is therefore an emphasis of this article to calculate all constants in painstaking detail.


Classical grid based methods in $d>1$ dimensions 
suffer from the \emph{curse of dimensionality}: the number of unknowns
$N$ required to achieve a prescribed error $\varepsilon$ grows exponentially in $d$  like
\[
N(\varepsilon) = \mathcal O (\varepsilon^{-d/p}),
\]
where $p$ is the order of the method.

The analysis in \citep{bungartz-griebel:04} shows that optimal approximation
 of sufficiently smooth functions, 
 for a given size of  a hierarchical tensor
product basis, is attained for so-called \emph{sparse grids}.
 Their relevance to the solution of PDEs was first revealed by 
\cite{zenger:90}, and subsequently error bounds for finite element methods for
elliptic problems were derived in detail by
\cite{bungartz:92, bungartz:98}.
More recently, optimal convergence rates of a sparse
wavelet method were obtained for parabolic equations, under
weak assumptions on the regularity of the initial condition,
using $hp$ discontinuous Galerkin time stepping in conjunction with the smoothing properties of parabolic equations \citep{petersdorff-schwab:02}.
What is more, adaptive sparse wavelet methods retain this optimal order for elliptic equations in the presence of corner singularities of the solution \citep{nitsche, stevenson:09, stevenson:10}.

In contrast to Galerkin-type methods, the combination technique -- first
introduced in \citep{griebel-schneider-zenger} -- decomposes the solution into
contributions from tensor product grids. This facilitates the discretisation of PDEs
has the practical
advantage that only numerical approximations on relatively small
conventional grids need to be computed. These can be obtained independently
and are superposed subsequently, which lends itself to very efficient parallel
implementations.
This concept has been successfully used in a number of applications, e.g.\
computational fluid dynamics \citep{griebel-thurner:95}, quantum mechanics
\citep{garcke-griebel:00} and computational finance \citep{reisinger:06}.

Theoretical results for this extrapolation scheme, however, which inevitably
rely on the expansion of the hierarchical surplus in terms of the grid sizes,
have so far only been obtained for simple model problems. 
\cite{bungartz-griebel-roeschke-zenger:94} obtain error bounds in terms of
the Fourier coefficients
for a central difference scheme for the two-dimensional Laplace equation.
\cite{pflaum:97} and \cite{pflaum-zhou:99} employ Sobolev space techniques for the
combination solution of a finite element method and prove
asymptotic errors of the form $h^2 |\log h|^{d-1}$ in the $L_2$-norm for general
linear elliptic equations in two dimensions, and first order convergence in the
$H^1$-norm. A similar result is also shown for the Poisson equation in higher
dimensions.
Common to these two approaches is the use of semi-discrete solutions to derive
expansions of the hierarchical surplus, namely Fourier
representations for semi-discrete solutions to the Laplace problem 
in \citep{bungartz-griebel-roeschke-zenger:94}, and
variational formulations with semi-discrete conforming
subspaces in \citep{pflaum:97, pflaum-zhou:99}.

We introduce a new framework to
derive error bounds for general difference schemes in arbitrary dimensions.
As in the aforementioned articles, again
semi-discrete problems will play an important role, now in the
somewhat different setting of an error correction scheme that
determines a suitable error expansion by a simple generic recursion.
This will be illustrated by the central difference
stencil for the Poisson problem and with upwinding for the advection
equation. For these two examples, sparse grid error bounds are then derived in
terms of the dimension of the problem and the smoothness of the solution. 
The results are formulated in spaces in which mixed
derivatives of sufficiently high order exist and vanish at the boundaries,
in order to avoid regularity problems in the corners.
We will discuss in the following section why this is not a problematic assumption for most practical applications where one would use sparse grids.


The rest of the paper is organised as follows.
Section \ref{subsec:combination} outlines the workings of the sparse grid combination technique,
the two main steps of the convergence analysis, and discusses smoothness requirements.
Section \ref{sec:error-expansions} develops a framework for the first 
problem-specific step of the analysis,
in which we derive a multivariate expansion of the form (\ref{asympt}) for the error
$u(\vf x_h)-\vf u_h$ between the numerical solution $\vf u_h$ and exact solution $u$
at the grid points $\vf x_h$. We derive exact coefficients for the Poisson problem and some extensions
to illustrate the application of the framework.
We then extend the grid solution to $[0,1]^d$ by
multilinear interpolation in Section \ref{sec:interpolation}, and show that
the expansion (\ref{asympt}) is maintained pointwise for the
error between the exact solution $u$ and the interpolated approximation
$\mathcal I \vf u_h$.
In the second step of the analysis, in Section \ref{sec:error-bounds}, these
error expansions on Cartesian grids are used
to estimate the error of the combined solution by combinatorial formulae.
Exact asymptotic expansions will be given.
In Section \ref{sec:numerical-results},
numerical results for carefully chosen model problems are discussed to illustrate the theory.
In particular, we highlight
the dependence of the error on the smoothness of the solution and the problem
dimension.
Section \ref{sec:discussion} discusses these results and points out future
directions, extensions to other problem classes and possible applications.

\section{The combination technique and main results}
\label{subsec:combination}

Consider the $d$-dimensional unit cube $I^d = [0,1]^d$ and a Cartesian grid with
mesh sizes $h_i = 2^{-l_i}$, corresponding to a level $l_i \in \N_0$ in
direction $i=1,\ldots,d$.

For a vector $\vf h = (h_i)_{1 \le i \le d}$ we denote by
$u_{\vf h}$ the approximation of a function on this grid with points
$\vf x_h = (i_j h_j)_{\scriptsize \begin{array}{l}0 \le i_j \le N_j, 1\le j\le
    d\end{array}}\!\!$, $N_j = 1/h_j = 2^{l_j}$.
Hereby $\vf u_h$ is the discrete vector of (approximated) values at the grid
points, which
is extended to $I^d$ by a suitable interpolation operator
$\mathcal I$ as $u_{\vf h} := \mathcal I \vf u_h$.
We will ultimately be concerned with a setting where $\vf u_h$ is the finite
difference solution to a scalar PDE, in which case $\vf u_h$ is different
from the vector obtained by evaluating the exact solution to the PDE, $u$, at
the grid points. Therefore we denote the latter by $u(\vf x_h)$.

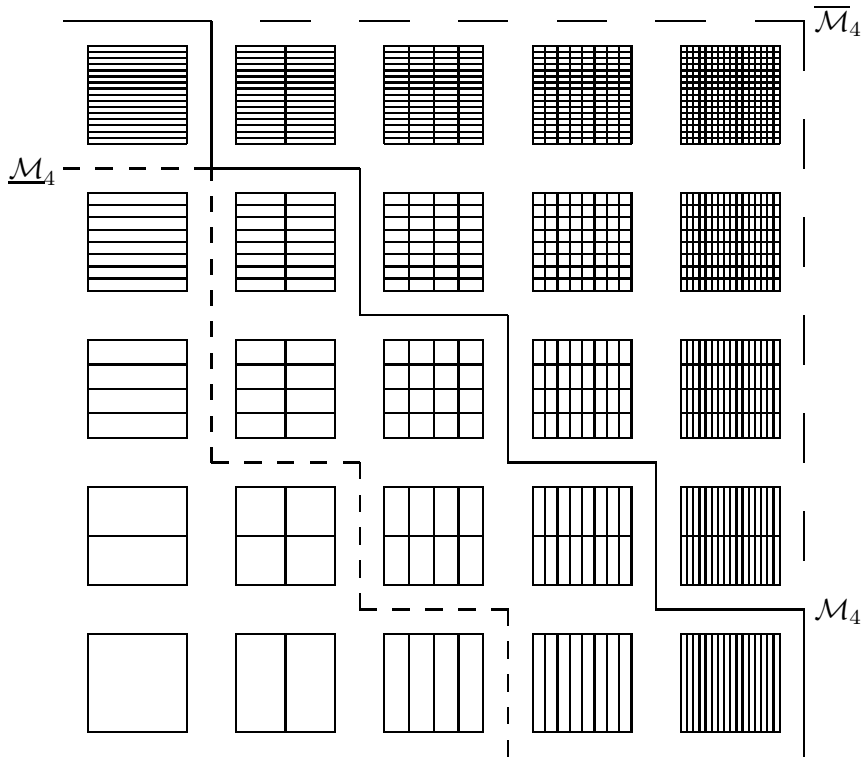
\begin{figure}
\centering
\setlength{\unitlength}{1.3 cm}
\newsavebox{\vertlinesO}
\savebox{\vertlinesO}(1,1)
{\begin{picture}(1,1)
\multiput(0,0)(1,0){2}{\line(0,1){1}}
\end{picture}}

\newsavebox{\vertlinesI}
\savebox{\vertlinesI}(1,1)
{\begin{picture}(1,1)
\multiput(0,0)(0.5,0){3}{\line(0,1){1}}
\end{picture}}

\newsavebox{\vertlinesII}
\savebox{\vertlinesII}(1,1)
{\begin{picture}(1,1)
\multiput(0,0)(0.25,0){5}{\line(0,1){1}}
\end{picture}}

\newsavebox{\vertlinesIII}
\savebox{\vertlinesIII}(1,1)
{\begin{picture}(1,1)
\multiput(0,0)(0.125,0){9}{\line(0,1){1}}
\end{picture}}

\newsavebox{\vertlinesIIII}
\savebox{\vertlinesIIII}(1,1)
{\begin{picture}(1,1)
\multiput(0,0)(0.0625,0){17}{\line(0,1){1}}
\end{picture}}

\newsavebox{\horilinesO}
\savebox{\horilinesO}(1,1)
{\begin{picture}(1,1)
\multiput(0,0)(0,1){2}{\line(1,0){1}}
\end{picture}}

\newsavebox{\horilinesI}
\savebox{\horilinesI}(1,1)
{\begin{picture}(1,1)
\multiput(0,0)(0,0.5){3}{\line(1,0){1}}
\end{picture}}

\newsavebox{\horilinesII}
\savebox{\horilinesII}(1,1)
{\begin{picture}(1,1)
\multiput(0,0)(0,0.25){5}{\line(1,0){1}}
\end{picture}}

\newsavebox{\horilinesIII}
\savebox{\horilinesIII}(1,1)
{\begin{picture}(1,1)
\multiput(0,0)(0,0.125){9}{\line(1,0){1}}
\end{picture}}

\newsavebox{\horilinesIIII}
\savebox{\horilinesIIII}(1,1)
{\begin{picture}(1,1)
\multiput(0,0)(0,0.0625){17}{\line(1,0){1}}
\end{picture}}

\begin{picture}(7,7)

\multiput(0,0)(0,1.5){5}{\usebox{\vertlinesO}}
\multiput(1.5,0)(0,1.5){5}{\usebox{\vertlinesI}}
\multiput(3,0)(0,1.5){5}{\usebox{\vertlinesII}}
\multiput(4.5,0)(0,1.5){5}{\usebox{\vertlinesIII}}
\multiput(6,0)(0,1.5){5}{\usebox{\vertlinesIIII}}

\multiput(0,0)(1.5,0){5}{\usebox{\horilinesO}}
\multiput(0,1.5)(1.5,0){5}{\usebox{\horilinesI}}
\multiput(0,3)(1.5,0){5}{\usebox{\horilinesII}}
\multiput(0,4.5)(1.5,0){5}{\usebox{\horilinesIII}}
\multiput(0,6)(1.5,0){5}{\usebox{\horilinesIIII}}


\multiput(7.25,-0.25)(0,1){8}{\line(0,1){0.5}}
\multiput(-0.25,7.25)(1,0){8}{\line(1,0){0.5}}


\put(7.25,-0.25){\line(0,1){1.5}}
\put(5.75,1.25){\line(0,1){1.5}}
\put(4.25,2.75){\line(0,1){1.5}}
\put(2.75,4.25){\line(0,1){1.5}}
\put(1.25,5.75){\line(0,1){1.5}}

\put(-0.25,7.25){\line(1,0){1.5}}
\put(1.25,5.75){\line(1,0){1.5}}
\put(2.75,4.25){\line(1,0){1.5}}
\put(4.25,2.75){\line(1,0){1.5}}
\put(5.75,1.25){\line(1,0){1.5}}


\multiput(4.25,-0.25)(0,0.33333333){5}{\line(0,1){0.16666}}
\multiput(2.75,1.25)(0,0.33333333){5}{\line(0,1){0.16666}}
\multiput(1.25,2.75)(0,0.32){11}{\line(0,1){0.153}}

\multiput(-0.25,5.75)(0.33333333,0){5}{\line(1,0){0.16666}}
\multiput(1.25,2.75)(0.33333333,0){5}{\line(1,0){0.16666}}
\multiput(2.75,1.25)(0.33333333,0){5}{\line(1,0){0.16666}}

\put(7.35,7.15){$\overline{\mathcal M}_4$}
\put(-0.8,5.65){$\underline{\mathcal M}_4$}
\put(7.35,1.15){${\mathcal M}_4$}

\end{picture}

\caption[Full grid, sparse grid, dimension adaptive sparse grid]{Full grid $\overline{\mathcal M}_4$, sparse grid $\mathcal M_4$,
  and a possible realisation of a dimension adaptive sparse grid
$\underline{\mathcal M}_4$.}
\label{gridtable}
\end{figure}

We now define the family $U$ of solutions corresponding to these grids (see
Fig.~\ref{gridtable}) by $U = (U(\vf i))_{\vf i \in \mathbb N_0^d}$ with
\[
U(\vf i) := u_{2^{-\vf i}},
\]
i.e.\ as the family of numerical approximations $u_{\vf h}$
on tensor product grids with $h_k = 2^{-i_k}$. Then the \emph{hierarchical
surplus} is the sequence
\begin{equation}
\label{surplus}
\del U := \del_1 \ldots \del_d U,
\end{equation}
where
\begin{eqnarray}
\label{delta}
\del_k U(\vf i)
:= \left\{\begin{array}{rl} U(\vf i)-U(\vf i -\vf e_k) & i_k>0, \\
U(\vf i) & i_k = 0,
\end{array}\right.
\end{eqnarray}
and
$\vf e_k$ 
the $k$-th unit
vector. The difference operators commute.

Given a sequence of index sets
$\mathcal{M}_n \subset \mathbb N_0^d$, the approximation on level $n$ is now defined as
\begin{equation}
\label{defsgsol}
u_n := \sum_{\vf i \in \mathcal{M}_n} \del U(\vf i).
\end{equation}
If this series converges in absolute terms, the error
\[
\|u_n - u\| = \Big\|\sum_{\vf i \notin \mathcal{M}_n} \del U(\vf i)\Big\|\le 
\sum_{\vf i \notin \mathcal{M}_n} \|\del U(\vf i)\|
\]
in a suitable norm $\|\cdot\|$
is therefore bounded by the surplus on finer grids which are not considered in the approximation
at level $n$.
The choice of an optimal refinement
strategy is thus determined by control of the surplus. If this is
done \emph{a posteriori}, dimension adaptive schemes are obtained \citep{Hegl2003,gerstner-griebel:03,griebelholtz}. This
requires that hierarchical surpluses on finer levels are estimated from
coarser levels. 
Conversely, an \emph{a priori} analysis requires analytic estimates of the
surplus in terms of the grid sizes $h_1, \ldots, h_d$ and is the
focus of this paper.

\cite{griebel-schneider-zenger} propose a splitting into lower-dimensional
contributions of the form
\begin{equation}
\label{asympt}
u-u_{\vf h} = \sum_{m=1}^{d} 
\sum_{
\tiny \begin{array}{c}
\{j_1,\ldots,j_m\} \\
\subset \{1,\ldots,d\}
\end{array}}
\gamma_{j_1,\ldots,j_m}
(\cdot;h_{j_1},\ldots, h_{j_m}) h_{j_1}^p\cdot \ldots \cdot h_{j_m}^p
\end{equation}
to analyse the combination technique. The `$\cdot$' stands for
the arguments of $u$. 
Before we
turn to the derivation of such expansions, we briefly outline the importance for the combination technique.

If we view the continuous function $u$ as a constant sequence over all refinement levels,
and therefore set $\delta u = 0$, then under assumption (\ref{asympt})
\begin{equation}
\label{order-delta}
\del U(\vf i) = \del \left(U(\vf i)- u \right) =
\mathcal{O}\left(2^{-p |\vf i|} \right),
\end{equation}
where $|\vf i| = |\vf i|_1 := \sum_{k=1}^d i_k$ and all lower order terms cancel
out through the difference operator.
Thus an asymptotically optimal choice for the index set for (\ref{defsgsol}) is
\[
\mathcal{M}_n = \{\vf i \in \N_0^d: |\vf i| \le n\}.
\]
The number of elements of $\mathcal{M}_n$ is then obtained by
\[
|\mathcal{M}_n/\mathcal{M}_{n-1}| = \bino{n+d-1}{d-1} \quad \Rightarrow \quad
|\mathcal{M}_n| = \sum_{l=0}^n \bino{l+d-1}{d-1} = \mathcal{O}(n^{d-1}),
\]
and the number of nodes in the grid corresponding to an index $\vf i$ is bounded by
$\prod_{k=1}^d (2^{i_k}+1) = \mathcal{O}\left(2^{|\vf i|}\right)$.
For this choice of grid we conclude from (\ref{order-delta}) that
\begin{eqnarray}
\label{esti}
\|u-u_n\| \le \sum_{\vf i \notin \mathcal M_n} \|\del U(\vf i)\| \le \sum_{|\vf
 i|>n} \mathcal O(2^{-p|\vf i|}) = \mathcal O(n^{d-1} 2^{-p n}).
\end{eqnarray}
The significance of the result (\ref{esti}) is that although the number of
degrees of freedom is only
\[
N_{dof} = \mathcal{O}\left(n^{d-1} 2^n \right),
\]
compared to $2^{d n}$ on the full grid, the error only deteriorates by a
factor of order $n^{d-1}$ compared to the full grid result $2^{-p n}$.

We illustrate the combination formula by a two-dimensional example:
\begin{eqnarray*}
u_0 &=& u(0,0) \\
u_1 &=& [u(1,0)-u(0,0)] + [u(0,1)-u(0,0)] + u_0 \\
&=& [u(1,0) + u(0,1)] - u(0,0) \\
u_2 &=& [u(2,0) - u(1,0)] + [u(1,1)-u(1,0)-u(0,1)+u(0,0)] + [u(0,2)-u(0,1)] +
u_1 \\
&=& [u(2,0) + u(1,1) + u(0,2)] - [u(1,0) + u(0,1)].
\end{eqnarray*}
The general structure in two dimensions is
\[
u_n = \sum_{l=0}^n U(l,n-l) - \sum_{l=0}^{n-1} U(l,n-1-l).
\]
In one dimension, the sparse grid is identical to the full grid.
In general, for $d\ge 1$ and $n>d-1$, the combined solution can be written as
\begin{eqnarray*}
u_n &=& [\dlt^{d-1} S](n), \\
S(n) &=& \sum_{|\vf i|=n} U(\vf i),
\end{eqnarray*}
where the difference operator $\dlt$, applied $d-1$ times to the
index $n$ of $S$, is the one-dimensional version of (\ref{surplus}), (\ref{delta}).
Evaluating the coefficients explicitly gives
\begin{equation}
  \label{combi}
  u_n = \sum_{l=n-d+1}^{n} \!\!\! a_{n-l} \; S(n),
\end{equation}
where
\begin{equation}
  \label{combi-coeff}
  a_i := (-1)^{d-1-i} \bino{d-1}{i}, \quad 0 \le i \le d-1.
\end{equation}

\label{subsec:def}

The error analysis falls into two main parts: If an expansion
(\ref{asympt}) can be shown for the problem at hand,
then error bounds of the form (\ref{esti})
follow by combinatorial arguments.

For the first problem-dependent part, we consider here linear PDEs on $I^d = [0,1]^d$ 
and denote by $\partial_i u$ the partial
derivative of $u$ with respect to $x_i$.
For a multi-index $\boldsymbol \alpha = (\alpha_1,\ldots,\alpha_d)$, let
$D^{\boldsymbol \alpha} u = \partial_1^{\alpha_1} \ldots
\partial_d^{\alpha_d} u$.
It will become clear later that
the appropriate function spaces are those with bounded mixed derivatives
in the supremum norm $\|\cdot\|_\infty$,
\begin{eqnarray}
X_{\boldsymbol \alpha}^d &:=& \left\{u: [0,1]^d \rightarrow \R : \;
D^{\boldsymbol \beta}u \in C_0(I^d) \; \forall\boldsymbol \beta \le \boldsymbol \alpha
\right \} \\
X_{\boldsymbol \alpha}^d(K) &:=& \left\{u \in X_{\boldsymbol \alpha}^d : \; 
\|D^{\boldsymbol \beta}u\|_\infty \le K \;
\forall\boldsymbol \beta \le \boldsymbol \alpha 
\right \} \\
\label{eqn:space}
X_k^d &:=& \left\{u: [0,1]^d \rightarrow \R : \;
D^{\boldsymbol \alpha}u \in C_0(I^d) \;
\forall\boldsymbol \alpha \in \{0,\ldots, k\}^d  
\right \} \\
X_k^d(K) &:=& \left\{u \in X_k^d: \;
\|D^{\boldsymbol \alpha}u\|_\infty \le K \;
\forall\boldsymbol \alpha \in \{0,\ldots, k\}^d
\right \}.
\end{eqnarray}

By $C_0(I^d)$ we denote the space of continuous functions which are zero
at the boundary, so we require for $X_k^d(K)$ that all mixed derivatives
up to order $k$ vanish at the boundary. 
We will now discuss the smoothness requirements,
and anticipate the following detailed results from the analysis in Section \ref{sec:error-bounds}.
\begin{theorem}
\label{theo-poisson-combined}
Let $u \in X_4^d$ be a solution of the Poisson problem and $u_n$ the
sparse grid solution with central differences on level $n$. Then
\begin{equation}
\|u-u_n\|_\infty \le c \cdot \sup_{|\boldsymbol \alpha|_\infty \le 4}
\|D^{\boldsymbol \alpha} u\|_\infty \cdot d \cdot
\left(\frac{5}{2}\right)^{d} \cdot (n+2(d-1))^{d-1} \cdot 4^{-n}
\label{final-poisson},
\end{equation}
where $c \le 121000$, $\boldsymbol \alpha = (4,\ldots,4)$, and
\begin{equation}
\label{poisson-final}
u-u_n = \tilde{c} \cdot \frac{d}{4^{3d}} \cdot
n^{d-1} \cdot 4^{-n} + \mathcal{O}\left(n^{d-2} 4^{-n} \right),
\end{equation}
for some $\tilde{c} \le c \|D^{\boldsymbol \alpha} u\|_\infty$.
\end{theorem}
\begin{theorem}
\label{theo-advection-combined}
Let $u \in X_2^d$ be a solution of the advection equation and $u_n$ the
sparse grid solution with upwinding and backward Euler timestepping on level $n$. Then
\begin{equation}
\|u-u_n\|_\infty \le c \cdot \sup_{|\boldsymbol \alpha|_\infty \le 2}
\|D^{\boldsymbol \alpha} u\|_\infty \cdot d \cdot
\left(\frac{3}{4}\right)^{d} \cdot (n+2(d-1))^{d-1} \cdot 2^{-n},
\label{final-advection}
\end{equation}
where $c \le 2$, $\boldsymbol \alpha = (2,\ldots,2)$, and
\begin{equation}
\label{advection-final}
u-u_n = \tilde{c} \cdot \frac{d}{4^d} \cdot
n^{d-1} \cdot 2^{-n} + \mathcal{O}\left(n^{d-2} 2^{-n} \right),
\end{equation}
for some $\tilde{c} \le c \|D^{\boldsymbol \alpha} u\|_\infty$.
\end{theorem}
The above theorems give error bounds of order $|\log h|^{d-1} h^p$,
and moreover make the asymptotic dependence of constants on the dimensionality and smoothness explicit. The analysis will show how these results generalise to a wider class of linear elliptic and parabolic PDEs, although an explicit computation of the constants is omitted in the general case.

Equations (\ref{poisson-final}) and (\ref{advection-final}) shows that for fixed $n$, the coefficients
of $n^{d-1} 2^{-p n}$ go to $0$ as $d\rightarrow \infty$, in line with observations in \citep{griebel06} and \citep{suli}.
This does not mean, however, that the approximation for a {given} level is better in higher
dimensions, due to the presence of the polynomial term in $n$. Moreover, this asymptotic approximation is obtained under
the assumption that $n$ is large compared to $d$, and this asymptotic range is
reached for increasing $n$ in higher dimensions, as
will be confirmed in the numerical results later.
This behaviour is also related to the fact that the minimum number of levels contained in
the combination solution grows linearly in $d$.

In the special case $\tilde{c} = 0$, the leading term proportional to $n^{d-1} 2^{-pn}$
vanishes. This is the case if the solution has a low
\emph{superposition dimension}, i.e.\ is a sum of functions that depend only
on a subset of the coordinates. 
In applications, this is often approximately the case, as can be seen by
asymptotic expansion \citep{reisinger:06}.
This effect can be exploited to construct generalised sparse grids to reduce the complexity to that of lower-dimensional sparse grids. 


We now return to the question of smoothness of solutions,
a key prerequisite for sparse grid approximation.
Due to the corners of the domain, the solution usually does not inherit the necessary degree of smoothness
from the data, e.g.\ for the Poisson problem, with homogeneous Dirichlet data on the $d$-cube, the solution for right hand-side $f=1$ does not have uniformly bounded mixed derivatives.
Rather, certain compatibility conditions for $f$ need to be satisfied.
For completeness,
conditions for sufficient regularity are derived in Appendix \ref{app:fourier},
and estimates for the derivatives of solutions are given in Appendix \ref{app:derivatives}.
This may seem to restrict the applicability of these results, and of sparse grids more generally, to
a small class of problems of limited practical value. We will now explain why this is not the case.

In most practically relevant high-dimensional applications, 
we are not interested in boundary value problems on the unit cube as such, but in equations on $\R^d$.
The localisation is necessary to make the problem computable on a grid, while a box shape of the computational domain is
amenable to tensor product grids. Boundary values are typically chosen by asymptotic analysis of the
unbounded problem to ensure that by making the box $B$ large enough, the difference between the original solution and the localised one can be made small. Now the unbounded solution
restricted to the box is smooth, and therefore the exact solution to the BVP with asymptotic boundary data, albeit usually not smooth itself, is a small perturbation of a smooth function.
Therefore, if the asymptotic boundary conditions are accurate to some $\epsilon$, and the
discretisation on each subgrid of the sparse grid is stable with respect to the boundary data,
which is normally given,
then the combined sparse grid error due to the boundary approximation is of order $\epsilon |\log h|^{d-1}$, where $|\log h|^{d-1}$ is the number of subgrids involved.
The error of the hypothetical discretisation of a BVP with exact boundary data on the cube is of order
$L^q h^p |\log h|^{d-1}$, as the continuous solution to this problem is smooth. 
The factor $L^q$, for some $q$ and with $L$ the size of the box, arises because of the transformation of the solution and its derivatives to the unit cube.
The sum of the localisation and discretisation errors is a bound for the discretisation error of the non-smooth localised solution. Indeed, we want to approximate the
unbounded solution, but we have assumed that the localisation error is bounded by $\epsilon$.
In many applications, asymptotic boundary values are known which converge exponentially to the true solution, as is seen e.g.\ for multi-dimensional Black-Scholes-type PDEs from
\citep{kangro}.
In these cases, one can therefore first pick the domain $B$ large enough and then $h$ to achieve a desired overall accuracy, without affecting the complexity order.

This justifies restricting the following error analysis to the setting where smooth solutions exist,
in particular we consider solutions which vanish sufficiently fast at the boundaries.



\section{Error expansion for finite difference schemes}
\label{sec:error-expansions}

The goal of this section is to establish error expansions of the form
\begin{equation}
\label{asympt2}
u(\vf x_h) - \vf u_{h} = \sum_{m=1}^{d} 
\sum_{
\tiny \begin{array}{c}
\{j_1,\ldots,j_m\} \\
\subset \{1,\ldots,d\}
\end{array}}
w_{j_1,\ldots,j_m}
(\vf x_h;h_{j_1},\ldots, h_{j_m}) h_{j_1}^p\cdot \ldots \cdot h_{j_m}^p,
\end{equation}
between the finite difference solution $\vf u_h$ and the exact solution $u$ of a PDE,
evaluated on a Cartesian grid $\vf x_h$ with grid sizes $h_1,\ldots,h_d$.

We have in mind linear (possibly degenerate) elliptic equations of the type
\[
A u = \sum_{i,j=1}^d a_{ij} \partial_i \partial_j u + \sum_{j=1}^d b_j
\partial_j u + c u = f,
\]
where $\sum_{i,j=1}^d n_i n_j a_{ij} \ge 0$,
and view parabolic equations as a degenerate case.

\subsection{Setup}

It seems instructive to first outline the main principles of the proof, which
are generic and can be applied to a large class of problems.
We write the discretised systems
in matrix notation
\begin{equation}
\label{discr-eqn}
\vf A_{h} \vf u_{h} = \vf f_{h},
\end{equation}
such that the truncation error can be written as $\vf A_h u(\vf x_h) - \vf f_h$.

In order to split up the truncation error into contributions from different dimensions, we define semi-discrete operators
$\vf A^{(i_1,\ldots,i_m)}_h$ in directions
$i_1,\ldots,i_m$, e.g.\ for the Laplace operator,
\begin{equation}
\label{semi-lapl}
\vf A^{(i_1,\ldots,i_m)}_h u = \sum_{i\in
  \{i_1,\ldots,i_m\}} \delta^+_{i,h_i} \delta^-_{i,h_i} u + \sum_{i\notin
  \{i_1,\ldots,i_m\}} \partial_i^2 u,
\end{equation}
where $\delta_{k,h_k}^+$ and $\delta_{k,h_k}^-$ are one-sided differences such that
$\delta^+_{i,h_i} \delta^-_{i,h_i} $ is a second difference.


In contrast to the fully discrete solutions $\vf u_h$ satisfying
(\ref{discr-eqn}),
this leads to a system of PDEs
\[
\vf A^{(i_1,\ldots,i_m)}_h \vf u^{(i_1,\ldots,i_m)}_h = \vf f^{(i_1,\ldots,i_m)}_h,
\]
with boundary conditions,
where $\vf u^{(i_1,\ldots,i_m)}_h$ are defined on hyper-planes
\[
\vf I^{(i_1,\ldots,i_m)} = \{x \in I^d: \; x_{i_k} \in \{j h_{i_k}: \; 0\le j
\le h_{i_k}^{-1}\}, 1\le k \le m \},
\]
and $\vf f^{(i_1,\ldots,i_m)}_h = \vf R^{(i_1,\ldots,i_m)}_h f$ is the
restriction of $f$ to $\vf I^{(i_1,\ldots,i_m)}$.
Let $C(\vf I^{(i_1,\ldots,i_m)})$
denote the space of continuous functions on these planes with maximum norm 
$\|\cdot\|_\infty$,
and derivatives are defined in directions along the planes.

We restrict the following sketch of the proof to two space dimensions in order to avoid
cumbersome notation.
We then extend this to higher dimensions and specify
the requirements for particular examples,
derive the coefficients in the expansion (\ref{asympt2}) and give sharp bounds.

\subsection{Outline of proof in two dimensions}

Starting point is a consistency assumption of order $p$ of the form
\begin{eqnarray}
\label{cons}
\vf A_h u - f = h_1^p \tau_1^{(0)}(\cdot,h_1) +
h_2^p \tau_2^{(0)}(\cdot,h_2) +
h_1^p h_2^p \tau_{1,2}^{(0)}(\cdot,h_1,h_2).
\end{eqnarray}
This will typically be straightforward to obtain
by Taylor expansion. In so doing, we assume that the solution is
sufficiently smooth, and
we make such an assumption throughout the following considerations.
Note that often, e.g.\ for the Poisson problem,
$\tau_{1,2}^{(0)}(\cdot,h_1,h_2) = 0$, but this term will be present, e.g., for the discretisation of mixed derivatives.

To derive the respective convergence order,
one would be tempted to write
\begin{eqnarray}
\label{naive}
u(\vf x_h) \!-\! \underbrace{\vf A_h^{-1} \vf f_h}_{\vf u_h} = 
h_1^p \vf A_h^{-1} \tau_1^{(0)}(\vf x_h,h_1) +
h_2^p \vf A_h^{-1} \tau_2^{(0)}(\vf x_h,h_2) +
h_1^p h_2^p \vf A_h^{-1} \tau_{1,2}^{(0)}(\vf x_h,h_1,h_2),
\end{eqnarray}
and deduce from the boundedness of $\vf A_h^{-1}$ convergence of the
scheme. In the present context, however, since $\vf A_h^{-1}$ depends on \emph{all} grid
sizes $h_1,\ldots,h_d$, it is not possible to derive an error expansion of the
form (\ref{asympt2}) in this way, unless $\vf A_h^{-1}$ is known explicitly in
terms of the grid sizes and an expansion of $\vf A_h^{-1}$ can be
derived. This is the principle behind the seminal work of
\cite{bungartz-griebel-roeschke-zenger:94}, where Fourier series of
continuous and discrete solutions to the Laplace equation are used. Note that
this is 
only possible in special cases where such representations are known.

This is where the concept of error correction comes into play. We 
determine the error terms by solving the auxiliary semi-discrete problems
\begin{eqnarray}
\label{semi-discr-sol}
\vf A_h^{(1)} w_{1}(\cdot,h_1) &=&
  \tau^{(0)}_{1}(\cdot;h_{1}) \\
\label{semi-discr-sol2}
\vf A_h^{(2)} w_{2}(\cdot,h_2) &=&
  \tau^{(0)}_{2}(\cdot;h_{2})
\end{eqnarray}
on the stacks of lines ${\vf I}_h^{(1)}$ and ${\vf I}_h^{(2)}$, respectively.
The crucial point is that $w_1(\cdot,h_1)$ and $w_2(\cdot,h_2)$ indeed only depend on $h_1$
and $h_2$ respectively,
because they are the solution of a semi-discretisation as in (\ref{semi-lapl})
and not of the fully discrete equation.

Under suitable regularity assumptions,
the first terms on the right-hand side of
\begin{eqnarray*}
\vf A_h \left(u - h_1^p w_1(\cdot;h_1) - h_2^p w_2(\cdot;h_2)\right) -
f &=& \\ && \hspace{- 6.7 cm}
h_1^p \left(\vf A_h^{(1)} w_1(\cdot;h_1) - \vf A_h w_1(\cdot;h_1)\right) +
h_2^p \left(\vf A_h^{(2)} w_2(\cdot;h_2) - \vf A_h w_2(\cdot;h_2)\right) +
h_1^p h_2^p \tau_{12}^{(0)}(\cdot;h_1,h_2)
\end{eqnarray*}
can be absorbed in the higher order terms by further expanding
\begin{eqnarray}
\label{semi-discr-exp}
\left[\vf A_h^{(1)} - \vf A_h \right] w_{1}(\cdot;h_1) &=&
h_2^p \, \sigma_{1;2}(\cdot;h_1;h_2), \\
\label{semi-discr-exp2}
\left[\vf A_h^{(2)} - \vf A_h \right] w_{2}(\cdot;h_2) &=&
h_1^p \, \sigma_{2;1}(\cdot;h_2;h_1),
\end{eqnarray}
using the truncation error of the semi-discrete problems. So if we define
\[
\tau^{(1)}_{1,2} = \tau^{(0)}_{1,2} +
\sigma_{1;2}(\cdot;h_1;h_2) + \sigma_{2;1}(\cdot;h_2;h_1),
\]
we get
\begin{eqnarray*}
\vf A_h \left(u - h_1^p \, w_1(\cdot;h_1) - h_2^p \, w_2(\cdot;h_2)\right) -
f =
h_1^p h_2^p \, \tau^{(1)}_{1,2}(\cdot;h_1,h_2).
\end{eqnarray*}
\emph{Now} we can conclude with a final step as in (\ref{naive}) that
\begin{eqnarray}
\label{result2d}
u(\vf x_h) - \vf u_h = h_1^p \, w_1(\vf x_h;h_1) + h_2^p\, w_2(\vf x_h;h_2)
+ h_1^p h_2^p\, w_{1,2}(\vf x_h;h_1,h_2),
\end{eqnarray}
where
\[
w_{1,2}(\vf x_h;h_{1},h_2) =  \vf A_h^{-1}
\tau^{(1)}_{1,2}(\vf x_h;h_{1},h_{2}).
\]

In higher dimensions, this procedure will be applied inductively.
The start of the induction, $m=1$, is always the consistency assumption
(\ref{cons}). The final step, $m=d+1$, is always essentially equivalent to
(\ref{result2d}).

The central proofs in this article go along these lines and it will become
clear how this framework can be applied to other settings.
The main gap in the above formal derivation concerns the existence of an expansion
(\ref{semi-discr-exp}), (\ref{semi-discr-exp2}), and the boundedness of its
coefficients, which relies on the smoothness of the solutions of
(\ref{semi-discr-sol}), (\ref{semi-discr-sol2}).
In the following two sections we detail the regularity requirements and give
explicit bounds on the coefficients in the expansion for the Poisson  problem and the
advection equation.

\subsection{Detailed analysis for the Poisson problem}

As an example, we consider central differences
for the Poisson problem 
\begin{eqnarray}
\label{poisson}
\Delta u = \sum_{i=1}^d \partial_i^2 u = f && \text{in } I^d, \\
u = 0 && \text{on } \partial I^d,
\end{eqnarray}
and assume sufficiently smooth compatible $f$ such that $u
\in X_4^d$.
The consistency order here is $p=2$.


Standard elliptic regularity results cannot be applied here 
as we need regularity of the solutions to
semi-discrete problems, 
with estimates independent of the grid sizes.
Therefore we formulate the following lemma.

\begin{lemma}
\label{lemma:stability}
Let $u$ be a solution of $\Delta u = f$ with homogeneous Dirichlet data and 
$\vf u_h^{(i_1,\ldots,i_m)}$ the solution of $\vf A_h^{(i_1,\ldots,i_m)} \vf
u_h^{(i_1,\ldots,i_m)} = \vf f_h^{(i_1,\ldots,i_m)}$ (with the definitions
from Section \ref{subsec:def}).
Then:
\begin{enumerate}
\item
$
\|u\|_\infty \le \frac{1}{8} \|f\|_\infty
$
\item
$
\|\vf u_h^{(i_1,\ldots,i_m)}\|_\infty \le \frac{1}{8} \|f\|_\infty
$
\end{enumerate}
\end{lemma}
\begin{proof}
\begin{enumerate}
\item
Let $|f|
< \bar{f}$
and $v := - \frac{1}{2} \bar{f} x_1 (1-x_1)$. Then
$
\Delta (u  + v) = f +  \bar{f} > 0
$
and hence $u + v < 0$ (maximum at the boundary), i.~e.\
$
u < -v \le \bar{f}/8.
$
Similarly $u\ge -\bar{f}/8$.
\item
Again a semi-discrete maximum principle holds and the analogous
result follows by considering
$
\vf A_h^{(i_1,\ldots,i_m)} v = \bar{f},
$
as central differences are exact for quadratic functions.
\end{enumerate}
\end{proof}

Similar estimates for the derivatives of the solution in terms of
the derivatives of $f$ can be obtained
by differentiating the equation, and using the fact that we are considering
function spaces for which derivatives of sufficiently high order vanish at
the boundaries.


For completeness, we now state the expansions (\ref{cons}) and
(\ref{semi-discr-exp}), (\ref{semi-discr-exp2}) in detail. Note that the truncation error
is defined continuously and not just at the grid points.

\begin{lemma}[Truncation error of difference stencil]
\label{lemma:consistency}
\begin{enumerate}
\item
Let $u \in X^d_{4}(K)$, then
\[
(A-\vf A_h) u = \sum_{k=1}^d h_k^2 \, \tau_k(\cdot;h_k)
\]
for some $\tau_k$ with
\begin{equation}
\label{estim-trunc}
\left\|D^{\boldsymbol \alpha} \tau_k \right\|_\infty \le \frac{1}{12} K
\end{equation}
for $\boldsymbol \alpha \in \{0,4\}^d$, $\alpha_k = 0$.
\item
Let $u \in X^d_{\boldsymbol \beta}(K)$ with $\beta_i = 4 \; \forall \; i
\notin \{i_1,\ldots,i_m\}$, then
\[
(\vf A_h^{(i_1,\ldots,i_m)}-\vf A_h) u \; = \sum_{k \notin \{i_1,\ldots,i_m\}}
h_k^2 \, \tau_k(\cdot;h_k),
\]
again with (\ref{estim-trunc}), but now for $\boldsymbol \alpha \in
\{0,4\}^d$, $\alpha_i = 0 \; \forall i \in \{i_1,\ldots,i_m\} \cup \{k\}$ and
$\alpha_i \le \beta_i$.
\end{enumerate}
\end{lemma}
\begin{proof}
Standard Taylor expansion in one variable.
\end{proof}


We now prove an error expansion for the finite difference solution at the
grid points.

\begin{theorem}
\label{theorem:expansion-on-grid}
Let $u \in X_4^d(K)$ be the solution of the Poisson equation and $\vf u_{h}$
the finite difference approximation on a grid $\vf x_h$. Then
\begin{equation}
\label{expansion-on-grid}
u(\vf x_h)-\vf u_{h} = \sum_{m=1}^{d} 
\sum_{
\scriptsize 
\begin{array}{c}
\{j_1,\ldots,j_m\} \\
\subset \{1,\ldots,d\}
\end{array}}
w_{j_1,\ldots,j_m}
(\vf x_h;h_{j_1},\ldots, h_{j_m}) h_{j_1}^2\cdot \ldots \cdot h_{j_m}^2,
\end{equation}
where $w_{j_1,\ldots,j_m} \in C\left(\vf I_h^{(j_1,\ldots,j_m)}\right)$ and
\[
\|w_{j_1,\ldots,j_m} (\cdot;h_{j_1},\ldots,
h_{j_m})\|_\infty \le K \frac{m!}{96^m}.
\]
\end{theorem}
\begin{proof}
We prove by induction for $1 \le m \le d$
\begin{eqnarray}
\nonumber
\vf A_h \left(u(\vf x_h) - \sum_{k=1}^{m-1}
\sum_{
\scriptsize 
\begin{array}{c}
\{i_1,\ldots,i_k\} \\
\subset \{1,\ldots,d\}
\end{array}}
h_{i_1}^2 \cdot \ldots \cdot h_{i_k}^2 w_{i_1,\ldots,i_k}(\vf
x_h;h_{i_1},\ldots,h_{i_k}) \right) -
  \vf f_h  &=& \\ && \hspace{-7.5 cm}
\sum_{
\scriptsize 
\begin{array}{c}
\{i_1,\ldots,i_m\} \\
\subset \{1,\ldots,d\}
\end{array}}
h_{i_1}^2 \cdot \ldots \cdot h_{i_m}^2 
\tau_{i_1,\ldots,i_m}(\vf x_h;h_{i_1},\ldots,h_{i_m}),
\label{induction}
\end{eqnarray}
where $w_{i_1,\ldots,i_k}(\cdot;h_{i_1},\ldots,h_{i_k})$ and
$\tau_{i_1,\ldots,i_m}(\cdot;h_{i_1},\ldots,h_{i_m})$
are functions defined on the hyper-planes ${\vf I}_h^{(j_1,\ldots,j_k)}$,
for which the estimates
\begin{eqnarray}
\label{tau-estimate}
\|D^{\boldsymbol \alpha} \tau_{i_1,\dots,i_m}(\cdot;h_{i_1},\ldots,h_{i_m})
\|_{\infty}
 &\le& m! \, 8^{-(m-1)} \, 12^{-m} \, K, \\
\|D^{\boldsymbol \alpha} w_{i_1,\dots,i_k}(\cdot;h_{i_1},\ldots,h_{i_k})\|_\infty
&\le& k! \, 8^{-k} \, 12^{-k} \, K \qquad \qquad \text{for } 1 \le k \le m-1
\label{w-estimate}
\end{eqnarray}
hold if $\boldsymbol \alpha \in \{0,4\}^d$ with $\alpha_{i_j} = 0$, $1 \le j
\le k$, i.e.\ along planes where derivatives are defined.

The case $m=1$,
\begin{eqnarray*}
\vf A_h u(\vf x_h) - \vf f_h = \sum_{k=1}^{d} h_{k}^2 \tau_{k}(\vf x_h, h_k),
\end{eqnarray*}
follows from Lemma \ref{lemma:consistency}, 1., with
$
\|D^{\boldsymbol \alpha} \tau_{k}\|_\infty \le \frac{1}{12} K
$
for $\boldsymbol \alpha \in \{0,4\}^d$ with $\alpha_k = 0$.
Now assume (\ref{induction}) holds for $m\ge 1$ with the bound (\ref{tau-estimate}).
Then the solution $w_{i_1,\ldots,i_m}$ of
\begin{equation}
\vf A_h^{(i_1,\ldots,i_m)} w_{i_1,\ldots,i_m} =
\tau_{i_1,\ldots,i_m}(\cdot;h_{i_1},\ldots,h_{i_m})
\end{equation}
satisfies, using Lemma
\ref{lemma:stability}, 2., and the comments thereafter,
\[
\|D^{\boldsymbol \alpha} w_{i_1,\dots,i_m}(\cdot;h_{i_1},\ldots,h_{i_m})\|_\infty
\le m! \, 8^{-m} 12^{-m} K,
\]
where
$\boldsymbol \alpha \in \{0,4\}^d$ and
$\alpha_{i_j} = 0$, $1 \le j \le k$.
Therefore, for $m<d$, from Lemma \ref{lemma:consistency}, 2., there exist
$\sigma_{i_1,\ldots,i_m;k}(\cdot;h_{i_1},\ldots,h_{i_m};h_{k})$ such that
\[
\vf A_h^{(i_1,\ldots,i_m)} w_{i_1,\ldots,i_m} - \vf A_h w_{i_1,\ldots,i_m} =
\sum_{k \notin \{i_1,\ldots,i_m\}} h_{k}^2
\sigma_{i_1,\ldots,i_m;k}(\cdot;h_{i_1},\ldots,h_{i_m};h_{k})
\]
and
\[
\|D^{\boldsymbol \alpha} \sigma_{i_1,\ldots,i_m;k}\|_{\infty} \le m! \, 8^{-m}
 12^{-(m+1)} K,
\]
for $\alpha_{i_j} = 0$, $1\le j\le m$ and $\alpha_k = 0$. For $m=d$, by construction
$\vf A_h^{(i_1,\ldots,i_m)}-\vf A_h = 0$.
So define
\[
\tau_{i_1, \ldots, i_{m+1}} (\cdot;h_{i_1},\ldots,h_{i_{m+1}})
:= \sum_{\scriptsize \begin{array}{c}j_1,\ldots,j_m,j \; {\rm s.t.}\\ 
\{j_1,\ldots,j_m\} \cup \{j\} = \{i_1, \ldots,i_{m+1} \}\end{array}}
\sigma_{j_1,\ldots,j_m;j}(\cdot;h_{j_1},\ldots,h_{j_m},h_j),
\]
and since the sum has $m+1$ terms,
\[
\|D^{\boldsymbol \alpha}
 \tau_{i_1,\ldots,i_m,i_{m+1}}\|_\infty
\le (m+1)! \, 8^{-m} 12^{-m-1} K,
\]
for $\boldsymbol \alpha \in \{0,4\}^d$ with $\alpha_{i_j} = 0$, $1 \le j \le
m+1$.
By induction, (\ref{induction}) follows for all $m$ and (\ref{expansion-on-grid}) 
is obtained directly by setting $m = d$.
\end{proof}

From the proof of Theorem \ref{theorem:expansion-on-grid} we see immediately
the following result:
\begin{corollary}
\label{corollary-smoothness}
The weights $w_{j_1,\ldots,j_m}(\vf x_h;h_{j_1},\ldots, h_{j_m})$ in 
(\ref{expansion-on-grid}) are the
restriction of functions $w_{j_1,\ldots,j_m}(\cdot;h_{j_1},\ldots, h_{j_m})$
defined on hyper-planes ${\vf I}_h^{(j_1,\ldots,j_k)}$, for which the bounds
\[
\|D^{\boldsymbol \alpha} w_{i_1,\dots,i_k}(\cdot;h_{i_1},\ldots,h_{i_k})\|_\infty
\le k! \, 96^{-k} \, K
\]
hold for $\boldsymbol \alpha \in \{0,4\}^d$,
$\alpha_i = 0 \; \forall i \in \{i_1,\ldots,i_k\}$.
\end{corollary}

\subsection{Advection equation and other extensions}

We now turn to the upwind discretisation of the advection equation
\begin{eqnarray}
\label{advection}
u_t + \sum_{i=1}^{d-1} \partial_i u = 0 && \quad \forall x \in I^{d-1},
\; t \in [0,1], \\
u(x,0) = u_0(x)  && \quad \forall x \in I^{d-1}, \\
u(x,t) = u_1(x,t) && \quad \forall x \in I^{d-1}, \; \exists i: \, x_i = 0, \; t
\in [0,1].
\label{end-advection}
\end{eqnarray}
For simplicity of notation, we take the velocity constant and equal to one in
each direction, but it will be obvious how the result generalises to the case
with variable velocity.
We can identify $t$ with $x_{d}$ in a combined space-time formulation to write
$\sum_{i=1}^d \partial_i u = 0$ with boundary condition
$u(x) = g(x) \; \forall x \text{ s.t. } \exists i : \, x_i = 0$.
Here $u_0$, $u_1$ (or $g$, respectively) are required to fulfill some compatibility conditions such that
$u \in X_2^d$ for the solution $u$ of (\ref{advection}).

Unconditional stability in the maximum norm is easy to
establish for the upwinding scheme
with left-sided
differences $\delta_{k,h_k}^-$ and implicit time-stepping.
For this instationary problem we will illustrate the
difference between a sparse grid in the ``space'' coordinates only and a ``space-time
sparse grid'' later.

Analogous results to Lemmas \ref{lemma:stability} and \ref{lemma:consistency}
can be derived and give the following expansion.

\begin{theorem}
\label{theorem:expansion-on-grid-adv}
Let $u \in X_2^d(K)$ be the solution of the advection equation and $\vf u_{h}$
the finite difference approximation. Then
\begin{equation}
\label{expansion-on-grid-adv}
u(\vf x_h)-\vf u_{h} = \sum_{m=1}^{d} 
\sum_{
\scriptsize 
\begin{array}{c}
\{j_1,\ldots,j_m\} \\
\subset \{1,\ldots,d\}
\end{array}}
w_{j_1,\ldots,j_m}
(\vf x_h;h_{j_1},\ldots, h_{j_m}) h_{j_1}\cdot \ldots \cdot h_{j_m},
\end{equation}
where $w_{j_1,\ldots,j_m} \in C\left(\vf I_h^{(j_1,\ldots,j_m)}\right)$ and
\[
\|w_{j_1,\ldots,j_m} (\cdot;h_{j_1},\ldots,
h_{j_m})\|_\infty \le K \frac{m!}{2^m}.
\]
\end{theorem}
\begin{proof}
First, derive an expression for the truncation error
(the equivalent of Lemma \ref{lemma:consistency}),
\[
(A-\vf A_h) u = \sum_{k=1}^d h_k \tau_k(\cdot;h_k),
\]
where
\[
\left\|D^{\boldsymbol \alpha} \tau_k \right\|_\infty \le \frac{1}{2} K
\]
for $\boldsymbol \alpha \in \{0,2\}^d$, $\alpha_k = 0$,
and a stability result (the equivalent of Lemma \ref{lemma:stability})
\[
\|\vf u_h^{(i_1,\ldots,i_m)}\|_\infty \le \max(\|u_0\|_\infty,\|u_1\|_\infty).
\]
The rest follows by the same steps as for the Poisson problem.
\end{proof}

We can combine the above results for the Poisson problem and the advection
equation to derive error formulae and estimates for an upwind discretisation of
advection-diffusion equations of the form
\[
\vf b \cdot \nabla u = \Delta u - c u.
\]
The truncation error is additive, regularity and stability apply similarly
under the assumptions made above. The derivation can therefore follow the same
steps.
Instationary problems
\[
u_t + \sum_{i=1}^d b_i \partial_i u = 
\sum_{i=1}^d a_i \partial_i^2 u - c u
\]
can be seen as a degenerate case with vanishing diffusion in one direction,
in which case upwinding in this coordinate is equivalent to fully
implicit time-stepping.
We can therefore construct a ``space-time'' sparse grid, which splits both
space and time in a hierarchical way. For other time discretisations, e.g.\
the explicit Euler scheme, which has a time step constraint
for stability in the maximum norm, time has to be treated separately:
a sparse grid is used for the space-like coordinates, while the time step has
to be chosen to satisfy the appropriate stability criterion on the highest
space refinement level.

For non-diagonal diffusion-tensors, it seems difficult to construct
finite difference schemes for which discrete maximum principles hold
on anisotropic grids.
In this case, the above analysis is not applicable. We discuss this in more
detail in Section \ref{sec:discussion}.

Broadly speaking,
all linear problems that admit smooth solutions and satisfy stability
properties that can be carried over to the discrete case fit into this
framework.

\section{Multilinear interpolation}
\label{sec:interpolation}

For the definition of an approximation on a sparse grid it is
necessary to extend the finite difference solution by interpolation. 
Interpolation on sparse grids clearly pre-dates any PDE connection \citep{smolyak},
and we only present results relevant to the further analysis.
We first show that the
interpolation error for a sufficiently smooth function by piecewise
multilinear splines on a Cartesian grid has an error
expansion of the form (\ref{error-expansion-approx}).
Subsequently we show that the difference between the numerical 
approximation of the PDE, i.e.\ the interpolated finite difference result, and the
exact solution still has such an expansion.

As a central ingredient, we first derive a particular partial Taylor
expansion. This resembles the semi-discretisations of the previous section.

\begin{lemma}[Expansion of functions with bounded mixed derivatives]
\label{lemma-expansion-high}
Let $u \in X_2^d$, $x=(x_1,\ldots,x_d)$, and introduce
$x^{(i_1,\ldots,i_k)} \in \R^d$
with $x^{(i_1,\ldots,i_k)}_j = x_j$ if $j\in \{i_1,\ldots i_k\}$, $0$
otherwise, and likewise for $s$, then
\begin{eqnarray*}
u(x) &=& u(0) \; + \; \sum_{i=1}^d x_i \, \partial_i u(x-x^{(i)})  \; +
\!\! \sum_{\scriptsize\begin{array}{c}i,j=1\\ i \neq j\end{array}}^d x_i x_j \,
\partial_i \partial_j u(x-x^{(i,j)})
+ \ldots  \\
&& \hspace{0.85 cm}
+\;\; x_1 \cdot \ldots \cdot x_d \, \partial_1 \ldots
\partial_d u(0) \; + \;
\sum_{i=1}^d \int_0^{x_i} (x_i-s_i)\, \partial_i^2 u(s^{(i)}) \ds_i +
\ldots \\
&& \hspace{0.85 cm}
+ \; \int_0^{x_1} \ldots \int_0^{x_d} 
(x_1-s_1)\cdot\ldots\cdot (x_d-s_d) \, \partial_1^2 \ldots \partial_d^2 u(s)
\ds_d \ldots \ds_1.
\end{eqnarray*}

%
%

\end{lemma}
\begin{proof}
See Appendix \ref{app:proofs3}.
\end{proof}

The obvious, but crucial point is that all terms that are linear in one or
more directions are interpolated exactly in these directions.

\begin{theorem}[Interpolation of functions with bounded mixed derivatives]
\label{lemma-approx}
Assume $u \in X_2^d$ and let $\mathcal{I} u(\vf x_h)$ be the multilinear
interpolating function on a grid $\vf x_h$. Then
\begin{equation}
\label{error-expansion-approx}
u(x)- (\mathcal{I} u(\vf x_h))(x) = \sum_{m=1}^{d} 
\sum_{
\tiny \begin{array}{c}
\{j_1,\ldots,j_m\} \\
\subset \{1,\ldots,d\}
\end{array}}
\alpha_{j_1,\ldots,j_m}
(x;h_{j_1},\ldots, h_{j_m}) h_{j_1}^2\cdot \ldots \cdot h_{j_m}^2,
\end{equation}
where
\[
\|\alpha_{j_1,\ldots,j_m} (\cdot;h_{j_1},\ldots, h_{j_m})\|_\infty \le
 \left(\frac{4}{27}\right)^m \left\|\partial_{j_1}^2 \ldots
 \partial_{j_m}^2 u \right\|_\infty.
\]
\end{theorem}
\begin{proof}
See Appendix \ref{app:proofs3}.
\end{proof}

\begin{remark}[Cubature]
From approximation results of this form, error expansions for cubature
formulae are obtained directly. Since the trapezoidal rule is exact on
piecewise
multilinear functions, the integration error is the integral of the
error terms over $[0,1]^d$.
For more detailed results see \cite{novak} or subsequent work, e.g.\
\cite{bungartz-griebel:04} and the references therein.
\end{remark}

To see that the error expansion 
(\ref{expansion-on-grid}) is preserved under interpolation of the
finite difference solution from the grid to $[0,1]^d$, we split
\begin{eqnarray}
\nonumber
u(x) - (\mathcal I \vf u_h) (x) 
&=& u(x) - (\mathcal I u(\vf x_h))(x) + (\mathcal I u(\vf x_h))(x) - 
(\mathcal I \vf u_h)(x) \\
&=& u(x) - (\mathcal I u(\vf x_h))(x) + (\mathcal I (u(\vf x_h) - \vf u_h))(x).
\label{error-splitting}
\end{eqnarray}
The first term $u(x) - (\mathcal I u(\vf x_h))(x)$ is the interpolation error of
the exact solution and is given by (\ref{error-expansion-approx}).
The second term, the interpolation of the discretisation error
\begin{equation}
\label{rep-expansion}
u(\vf x_h) - \vf u_h =
\sum_{k=1}^{d} \sum_{\tiny \begin{array}{c}
\{j_1,\ldots,j_m\} \\
\subset \{1,\ldots,d\}
\end{array}} h_{j_1}^2 \cdot \ldots \cdot h_{j_k}^2
w_{j_1,\ldots,j_k}(\vf x_h; h_{j_1},\ldots,h_{j_k})
\end{equation}
on the grid $\vf x_h$, is problem dependent
(see Theorem \ref{theorem:expansion-on-grid}) and needs
to be evaluated separately. It is straightforward to see that the linear
interpolant of the error is bounded by the error at the grid points, but it is
essential, and more involved, to derive the exact form of this expansion.

We show this again for the Poisson problem first.

\begin{theorem}
\label{theorem:error-expansion-total}
Let $u \in X_4^d(K)$ be the solution to the Poisson problem and $\vf u_h$ the
numerical solution with central differences, $\mathcal{I} \vf u_h$ its
multilinear interpolation. Then
\begin{equation}
\label{error-expansion-total}
u - \mathcal{I} \vf u_h = \sum_{m=1}^{d}
\sum_{
\tiny \begin{array}{c}
\{j_1,\ldots,j_m\} \\
\subset \{1,\ldots,d\}
\end{array}}
v_{j_1,\ldots,j_m}
(\cdot;h_{j_1},\ldots, h_{j_m}) h_{j_1}^2\cdot \ldots \cdot h_{j_m}^2,
\end{equation}
where
\[
\|v_{j_1,\ldots,j_m}
(\cdot;h_{j_1},\ldots, h_{j_m})\|_\infty \le C \cdot K
\cdot \frac{m!}{96^m},
\]
and $C < 150188$ depends on neither the dimension nor the data.
\end{theorem}
\begin{proof}
See Appendix \ref{app:proofs3}.
\end{proof}

Similarly, one gets for the advection equation the following result.
\begin{theorem}
\label{theorem:error-expansion-total-upw}
Let $u \in X_2^d(K)$ be the solution of the advection equation and $\vf u_h$ the
numerical solution with upwinding and implicit Euler time-stepping. Then
\begin{equation}
\label{error-expansion-total-upw}
u - \mathcal{I} \vf u_h = \sum_{m=1}^{d}
\sum_{
\tiny \begin{array}{c}
\{j_1,\ldots,j_m\} \\
\subset \{1,\ldots,d\}
\end{array}}
v_{j_1,\ldots,j_m}
(\cdot;h_{j_1},\ldots, h_{j_m}) h_{j_1} \cdot \ldots \cdot h_{j_m},
\end{equation}
where
\[
\|v_{j_1,\ldots,j_m}
(\cdot;h_{j_1},\ldots, h_{j_m})\|_\infty \le C \cdot K
\cdot \frac{m!}{2^m},
\]
and $C < 3/2$ depends on neither the dimension nor the data.
\end{theorem}
\begin{proof}
See Appendix \ref{app:proofs3}.
\end{proof}

\section{Combined error bounds and asymptotics}
\label{sec:error-bounds}

We study now in detail, by means
of combinatorial relations, the extrapolation effect that the combination formula has on error
terms. This leads to error bounds for the combination
solution of the order seen in (\ref{esti}).
Starting point is the pointwise expansion of the error on tensor product grids
with mesh sizes $\vf h = (h_1,\ldots,h_d)$,
\begin{equation}
\label{satz-asympt}
u-u_{\vf h} = \sum_{m=1}^d 
\sum_{\scriptsize \begin{array}{c}\{j_1,\ldots, j_m\} \\
\subset \{1,\ldots,d\} \end{array}}
v_{j_1,\ldots,j_m}
(\cdot;h_{j_1},\ldots, h_{j_m}) h_{j_1}^p\cdot \ldots \cdot h_{j_m}^p,
\end{equation}
where
\begin{equation}
\label{Kest}
|v_{j_1,\ldots,j_m}| \le K \qquad \forall 1\le m\le d \;\; \forall
\{j_1,\ldots,j_m\} \subset \{1,\ldots,d\},
\end{equation}
as shown for the examples of the preceding sections.
This second step of the analysis, however, is independent of the problem,
given that (\ref{satz-asympt}) holds.

\label{subsec:error-bounds}

\cite{griebel-schneider-zenger}  derive from
(\ref{satz-asympt}) and (\ref{Kest}),
for $d=2, 3$ and $p=2$, the bounds
 \begin{equation}
 \label{griebel-2d}
 |u-u_n| \le K 2^{-2n} \left(1 + \frac{5}{4} n \right)
 \end{equation}
and
 \begin{equation}
 \label{griebel-3d}
 |u-u_n| \le K 2^{-2n} \left(1 + \frac{65}{32} n
 + \frac{25}{32} n^2 \right),
 \end{equation}
respectively.
This section is devoted to a generalisation of (\ref{griebel-2d}) and
(\ref{griebel-3d}) to error bounds of the form
\[
|u-u_n| \le K c(d,p) n^{d-1} 2^{-pn},
\]
for arbitrary dimension and arbitrary order.
In addition, the asymptotic limit
\[
\lim_{n\rightarrow \infty} \frac{|u-u_n|}{n^{d-1} 2^{-pn}}
\]
will be given.
We will use the representation of the combined solution
\begin{eqnarray}
\label{combi-sol}
u_n &=& [\dlt^{d-1} S](n), \\
S(n) &=& \sum_{|\vf i|=n} U(\vf i),
\end{eqnarray}
as seen in Section \ref{subsec:combination}, in terms of an iterated application
of the one-dimensional difference operator $\dlt$.
The proof requires a few combinatorial identities. The following
formula for iterated differences of products --- a discrete version of the
product rule for differentiation --- will be useful.
The straightforward proof is omitted here, but see \citep{reisinger:04}.
 \begin{proposition}
 \label{prod-diff}
Let $f, g \in \R^{\N_0}$. Then
 \begin{equation}
 \dlt^k \left(f g\right) = \sum_{j=0}^k \bino{k}{j}
 \dlt^{k-j} f(\cdot+j) \dlt^j g \quad \forall k\in \N_0.
 \end{equation}
 \end{proposition}

Since the number of grids on level $n$ involved in the combination
solution in dimension $d$ is given by
\begin{equation}
\label{Ndef}
N(n,d) = \bino{n+d-1}{d-1},
\end{equation}
the following Lemma \ref{konsistenz}
states a necessary condition for the
consistency of the combination technique (i.e.~the sum of coefficients of all
grids is $1$).
 \begin{lemma}[Consistency]
 \label{konsistenz}
With $N$ from (\ref{Ndef}) $\forall d\in \N$, $\forall n\in \N, n\ge d-1$
 \[
\dlt^{d-1} N(n,d) = 1.
 \]
 \end{lemma}
 \begin{proof}
By induction one proves for $0 \le k \le d-1$
 \[
 \dlt^k N(n,d) = \bino{n+d-1}{d-k-1}.
 \]
 \end{proof}

The following formula is key to the proof of the main result in this
section, Theorem \ref{zweiteOrdnung}.
\begin{lemma}[Error representation formula]
\label{err-rep}
Let $m,d \ge 1$, 
$v: \R_+^m \rightarrow \R$ and for $n
\in \N_0$
\[
F(n) := \sum_{\scriptsize \begin{array}{c} \vf i \in \N_0^d \\ |\vf{i}|=n
  \end{array}}
 v(2^{-i_1},\ldots,2^{-i_m})
2^{-p i_1} \cdot \ldots \cdot 2^{-p i_m}.
\]
Then, for $d \in \N$,
\[
\dlt^{d-1} F(n) =
2^{-p(n+d-1)} \sum_{i=0}^{m-1} s_{n+d-i-1} \bino{m-1}{i}
(-2)^{pi},
\]
where
$
s_l := \sum_{\scriptsize \begin{array}{c} \vf i \in \N_0^d \\ |\vf{i}|=l
  \end{array}} v(2^{-i_1},\ldots,2^{-i_m}).
$
\end{lemma}
\begin{proof} See Appendix \ref{app:proofs4}. \end{proof}
Consider in Lemma \ref{err-rep} $n$ as the sparse grid level and
$s_l 2^{-pl}$ as the
order $p$ error on level $l$, with some coefficient $s_l$ that is collected
from a number of grids given by the
binomial term. Then Lemma \ref{err-rep} says that only the highest order terms
$2^{-p(n+d-1)}$ are left in the sparse grid solution, whereas
the combination formula cancels out all lower
order terms that come from the larger mesh sizes on the anisotropic grids.
This is a more quantitative version of (\ref{order-delta}).
After these preparations, the error terms can be estimated conveniently.

\begin{theorem}[Error bounds]
\label{zweiteOrdnung}
Assume for all $u_{\vf{h}}$ a pointwise error expansion of the form
(\ref{satz-asympt}) with (\ref{Kest}).
Then the combination solution (\ref{combi-sol}) fulfills the error estimate
\begin{equation}
\label{abschZweite}
|u-u_n| \le \frac{2 K}{(d-1)!} \left(\frac{2^p+1}{2^{p-1}}\right)^{d-1}
 (n+2(d-1))^{d-1} 2^{-pn}.
\end{equation}
\end{theorem}
\begin{proof}
See Appendix \ref{app:proofs4}.
\end{proof}

Let us study equation (\ref{abschZweite}) for $p=2$ and small $d$. For $d=1$, where
the sparse grid is identical to the full grid, the estimate reduces to
$|u-u_n| \le 2 K 4^{-n}$. The substitution $2^d-1$ by $2^d$ in the end of the
above proof explains the unnecessary factor 2.
For $d=2$, the leading term is given by (\ref{dritteUngl}), such that
in the highest power of $n$ 
\[
|u-u_n| \sim K \frac{5}{4} n 4^{-n},
\]
in accordance with (\ref{griebel-2d}).
Similarly, in three dimensions one gets
\[
|u-u_n| \sim K \frac{1}{2} \left(\frac{5}{4}\right)^2  4^{-n} =
K \frac{25}{32}  n^2 4^{-n}
\]
and recovers  (\ref{griebel-3d}).
Lower order terms differ due to the numbering of the grids and in fact they were
not estimated separately from (\ref{zweiteUngl}) onwards.

More interesting, however, is the dependence of (\ref{abschZweite}) on the
dimension $d$ for large $d$.
Keeping $n$ fixed, one sees from Stirling's formula,
\[
k! \sim \sqrt{2 \pi k} \left(\frac{k}{e}\right)^k\! ,
\]
that the factor depending explicitly on $d$ grows asymptotically like
$(5 \e)^d/\sqrt{d}$.
It is an interesting and practically very relevant question, whether the bounds
in (\ref{abschZweite}) are sharp asymptotically or the
exponentially growing constants can be omitted by
subtle treatment of the error terms.



\begin{corollary}[to Theorem \ref{zweiteOrdnung}]
\label{supiAbsch}
Under the assumptions of Theorem \ref{zweiteOrdnung} the sharper bound
\[
|u-u_n| \le 2 K \left(\frac{2^p+1}{2^{p-1}}\right)^{d-1} \!\!
\left(1+ (n+d-1) \frac{1+\ln(d-1)}{d-1}\right)^{d-1} \! 2^{-p n}
\]
holds for $d\ge 2$.
\end{corollary}
\begin{proof} See Appendix \ref{app:proofs4}. \end{proof}

\begin{corollary}[to Theorem \ref{zweiteOrdnung}]
\label{exactAsympt}
If additionally $v_{1,\ldots,d}$ in (\ref{satz-asympt}) is continuous in
$\vf 0 \in \R^d$ with
$\bar{v} := v_{1,\ldots,d}(\cdot;0,\ldots,0) \ne 0$, then the asymptotic
behaviour is
\begin{equation}
\label{traum}
u-u_n = \bar{v} \left(\frac{2^p-1}{2^p}\right)^{d-1}
 \frac{n^{d-1}}{(d-1)!} 2^{-p n} + \mathcal{O}\left(n^{d-2} 2^{-pn}\right).
\end{equation}
\end{corollary}

\begin{proof} See Appendix \ref{app:proofs4}.
\end{proof}


Finally we can collect the results 
from
Theorems 
\ref{zweiteOrdnung} and
\ref{theorem:error-expansion-total},
and
Theorems \ref{zweiteOrdnung} and \ref{theorem:error-expansion-total-upw}, respectively,
 to obtain Theorems \ref{theo-poisson-combined} and \ref{theo-advection-combined}.

\section{Numerical results}
\label{sec:numerical-results}
We illustrate the theoretical findings by numerical experiments, and pay
particular attention to the asymptotic convergence order, the dependence of
the error on the dimension and on the smoothness of the solution, as reflected in
(\ref{final-poisson}) and (\ref{final-advection}).

\subsection{Elliptic  problems}

Consider central differences for
\begin{eqnarray}
\label{eqn-poisson}
\Delta u &=& f \quad \text{ in } [0,1]^d, \\
u &=& g \quad \text{ on } \partial [0,1]^d.
\nonumber
\end{eqnarray}
We choose the data such that the solution is
\begin{equation}
\label{sol-poisson}
u(x) = \exp \left(- \frac{1}{2} \sum_{i=1}^d \lambda_i (x_i-p_i)^2 \right),
\end{equation}
where $\lambda_i \ge 0$ and $p \in [0,1]^d$, that is
$f(x) = \sum_{i=1}^d \lambda_i \left(-1 +\lambda_i y_i^2 \right) \cdot u(x).
$
This allows us to control all derivatives in the light of (\ref{final-poisson}).
In particular,
\[
\|D^{2 \boldsymbol \alpha} u\|_\infty = |D^{2 \boldsymbol \alpha} u(p)| =
\prod_{i=1}^d \lambda_i^{\alpha_i}.
\]
We choose $p_1 = 0.22081976$, $p_2 = 0.29072005$, $p_3 = 0.28051979$,
$p_4 = 0.27032006$, $p_5 = 0.24122005$, $p_6 = 0.17071947$,
$p_7 = 0.10101947$, $p_8 = 0.09021981$ to avoid symmetry effects.



We start by considering (\ref{sol-poisson}) with
$d \ge 1$ and $\lambda_i = 1$, $i = 1,\ldots,d$. Then, from above,
$\max_{|\boldsymbol \alpha|_\infty \le 4} \|D^{\boldsymbol \alpha} u\|_\infty$ = 1
for all $d$.

\begin{center}
\begin{figure}[ht]
\psfrag{n}[l][r][1.0]{$n$}
\psfrag{err}[l][u][1.0]{$\log_{2} \epsilon_n$}
\begin{center}
\includegraphics[width= 0.48 \columnwidth]{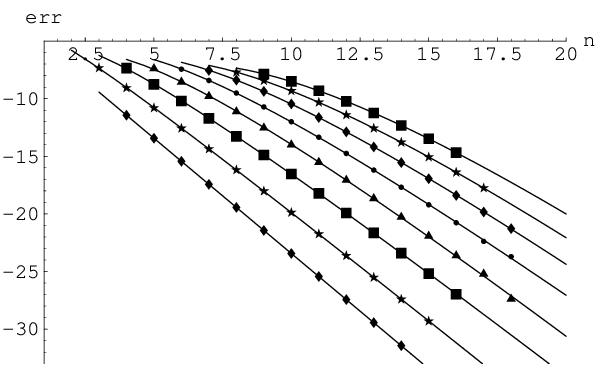} \hfill
\psfrag{err}[l][u][1.0]{$\log_{2} \epsilon_N$}
\psfrag{n}[l][r][1.0]{$\log_2 N$}
\includegraphics[width= 0.48 \columnwidth]{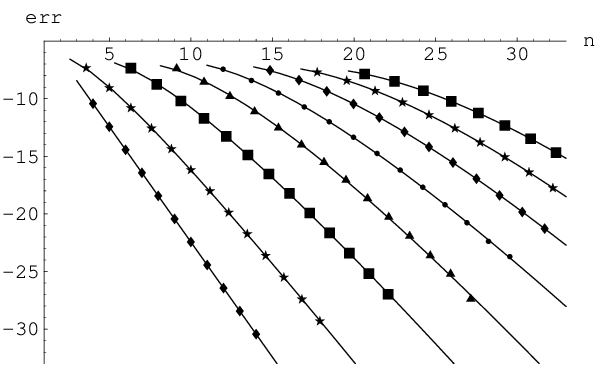}
\end{center}
\caption[Pointwise sparse grid error in dimensions 1 up to 8]{Pointwise error $\epsilon_n$ vs.\ grid level $n$ (left) and vs.\ the number of
grid points $N$ (right) for a sparse grid in dimensions
1 up to 8 (from left to right). The continuous curves are the asymptotes as
in (\ref{asymp-curve}), fitted to the data -- the estimated parameters are given in
Table \ref{tab:coeffs}.}
\label{fig:errorDim}
\end{figure}
\end{center}

Figure \ref{fig:errorDim} shows a logarithmic plot of the error
\begin{equation}
\label{pt-error1}
\epsilon_n := |u(x_*) - u_n(x_*)|
\end{equation}
for the sparse grid solution $u_n$ at level $n$, evaluated at a fixed point $x_*$.
Here $x_* = (1/2,\ldots,1/2)$, the centre point.

\begin{center}
\begin{figure}[ht]
\psfrag{n}[l][r][1.0]{$n$}
\psfrag{err}[l][u][1.0]{$\log_{2} \epsilon_n$}
\begin{center}
\includegraphics[width= 0.48 \columnwidth]{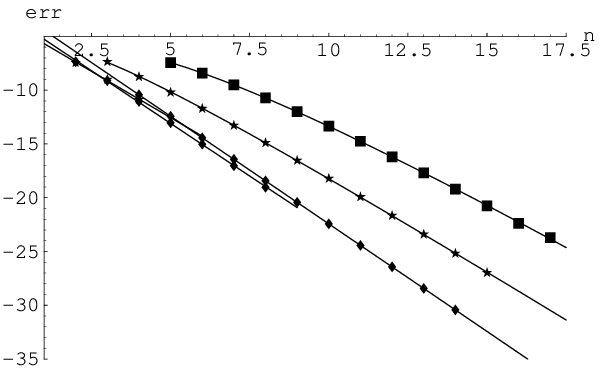} \hfill
\psfrag{err}[l][u][1.0]{$\log_{2} \epsilon_N$}
\psfrag{n}[l][r][1.0]{$\log_2 N$}
\includegraphics[width= 0.48 \columnwidth]{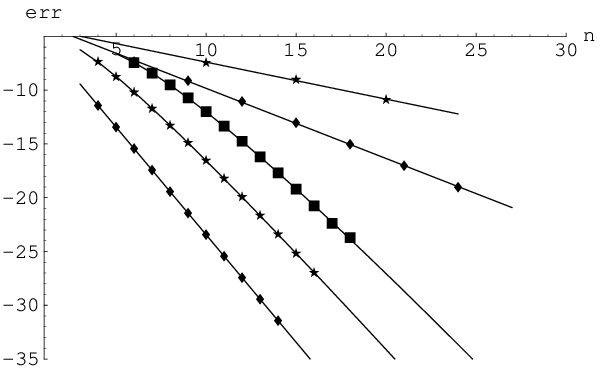}
\end{center}
\caption[Comparison sparse and full grid]{A similar plot to Fig.~\ref{fig:errorDim}, but comparing the sparse
grid in dimensions 1, 3 and 5 (upper graphs in the left plot and lower graphs on the
right) to the full grid. The full grid error on a given level does not
depend significantly on the dimension (three lines on top of each other), whereas
the sparse grid requires more refinements to achieve the same accuracy.
The right plot, however, shows the superior complexity in terms of error reduction per
unknown.}
\label{fig:errorDimSF}
\end{figure}
\end{center}

From the result (\ref{poisson-final}) we know that
\[
\epsilon_n = \mathcal{O}(n^{q} 2^{-p n})
\]
with $q = d-1$ and $p=2$. We therefore fit the function
\begin{equation}
\label{asymp-curve}
l(n,p,q,r) = -r + q \log_2 n - p n
\end{equation}
to $\log_2 \epsilon_n$ and determine $r$, $p$ and $q$ by least squares.
The results shown in Table \ref{tab:coeffs} correspond very well to the theoretical
values. It is clearly difficult to estimate the exponent of the logarithm to good accuracy.
\begin{table}
\centering
\begin{tabular}{|l|r|r|r|r|r|r|r|r|}\hline
$d$ & 1 & 2 & 3 & 4 & 5 & 6 & 7 & 8 \\ \hline
$p$ & 2.000 (2) & 1.938 (2) & 1.905 (2) & 1.970 (2) & 1.871 (2)
& 1.901 (2) & 1.944 (2) & 1.982 (2) \\ \hline
$q$ & -0.001 (0) & 0.478 (1) & 1.44 (2) & 2.86 (3) & 3.42 (4) & 4.75 (5) &
6.24 (6) & 7.76 (7) \\ \hline
$r$ & 11.43 & 3.97 & 3.97 & 5.31 & 6.04 & 8.42 & 11.61 & 15.28 \\ \hline
\end{tabular}
\caption[Fitted error coefficients]{Coefficients as in (\ref{asymp-curve}), fitted to the computed errors
by regression. In brackets see the values from the theory. The corresponding
curves are plotted with the data in Fig.~\ref{fig:errorDim}.}
\label{tab:coeffs}
\end{table}

We perform the same exercise for the error versus the number of grid points
$N$ and fit
\begin{equation}
\label{asymp-curve2}
l(N,\tilde{p},\tilde{q},\tilde{r}) = -\tilde{r} + \tilde{q} \log_2 N - \tilde{p} N
\end{equation}
to the data depicted on the right of Fig.~\ref{fig:errorDim} and
Fig.~\ref{fig:errorDimSF}. The coefficients for the sparse grid are given in
Table \ref{tab:coeffs2} and show that because the asymptotic complexity of the
sparse grid is independent of the dimension, up to logarithmic factors, the order
remains approximately 2.

\begin{table}
\centering
\begin{tabular}{|l|r|r|r|r|r|r|r|r|}\hline
$d$ & 1 & 2 & 3 & 4 & 5 & 6 & 7 & 8 \\ \hline
$\tilde{p}$ & 2.000 & -1.889 & 1.815 & 1.846 & 1.715 & 1.702 & 1.689 & 1.668 \\ \hline
$\tilde{q}$ & -0.002 & 2.23 & 5.08 & 8.74 & 10.6 & 13.7 & 16.7 & 19.6 \\ \hline
$\tilde{r}$ & 2.43 & 4.75 & 9.48 & 18.5 & 24.9 & 35.5 & 47.1 & 59.2 \\ \hline
\end{tabular}
\caption[Fitted complexity coefficients]{Coefficients as in (\ref{asymp-curve2}), fitted to the computed errors
by regression.}
\label{tab:coeffs2}
\end{table}

For comparison, the lines for the full grid in dimensions 1, 3 and 5 have slope $-2$
(identical to the sparse grid), $- 0.65$ (asymptotically $-2/3$) and $-0.34$
(asymptotically $-2/5$).
This reflects the curse of dimensionality, $\epsilon = \mathcal{O}(N^{-2/d})$.

\begin{center}
\begin{figure}[ht]
\psfrag{eps}[r][l][1.0]{$\log_2 \epsilon_n$}
\psfrag{2}[c][c][0.8]{$2$}
\psfrag{4}[c][c][0.8]{$4$}
\psfrag{6}[c][c][0.8]{$6$}
\psfrag{8}[c][c][0.8]{$8$}
\psfrag{10}[c][c][0.8]{$10$}
\psfrag{12}[c][c][0.8]{$12$}
\psfrag{14}[c][c][0.8]{$14$}
\psfrag{-5}[r][r][0.8]{$-5$}
\psfrag{-10}[r][r][0.8]{$-10$}
\psfrag{-15}[r][r][0.8]{$-15$}
\psfrag{-20}[r][r][0.8]{$-20$}
\psfrag{n}[r][l][1.0]{$n$}
\psfrag{plot1}[l][l][0.6]{$\lambda_1 = 1$, $\lambda_2 = 1$, $\lambda_3 = 1$}
\psfrag{plot2}[l][l][0.6]{$\lambda_1 = 100$, $\lambda_2 = 1$, $\lambda_3 = 1$}
\psfrag{plot3}[l][l][0.6]{$\lambda_1 = 10$, $\lambda_2 = 20$, $\lambda_3 = 50$}
\psfrag{plot4}[l][l][0.6]{$\lambda_1 = 100$, $\lambda_2 = 100$, $\lambda_3 = 1$}
\psfrag{plot5}[l][l][0.6]{$\lambda_1 = 10000$, $\lambda_2 = 1$, $\lambda_3 = 1$}
\begin{center}
\includegraphics[width= \columnwidth]{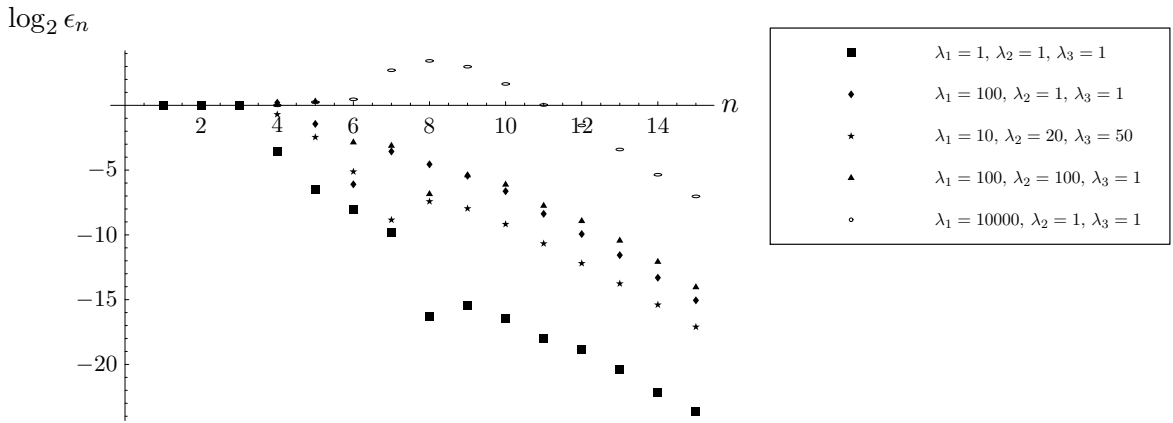}
\end{center}
\caption[Pointwise sparse grid error for anisotropic solutions]{Pointwise sparse grid error $\epsilon_n$ as in (\ref{pt-error1})
on level $n$ for equation (\ref{eqn-poisson}) with solution (\ref{sol-poisson})
for different values for $\lambda_1$, $\lambda_2$ and $\lambda_3$.}
\label{fig:lambdaPlotsConv}
\end{figure}
\end{center}

It remains to study the effect of smoothness and (an-)isotropy on the
convergence. We consider the three-dimensional case and vary $\lambda_1$, $\lambda_2$
and $\lambda_3$. The leading error term is proportional to
\[
\max_{x \in [0,1]^d}\frac{\partial^{12} u}
{\partial x_1^4 \partial x_2^4 \partial x_3^4}(x) = \lambda_1^2
\lambda_2^2 \lambda_3^2.
\]
The evaluation point $x_*$ is chosen equal to the point $p$.
To illustrate the scales involved, we plot $u$ for a selection of parameters
with fixed $\lambda_3 = 1$ and for a cross-section $x_3 = p_3$.
\begin{center}
\begin{figure}[ht]
\psfrag{u}[r][l][1.0]{$u$}
\psfrag{x1}[c][c][1.0]{$x_1$}
\psfrag{x2}[l][u][1.0]{$x_2$}
\begin{center}
\includegraphics[width= 0.3 \columnwidth]{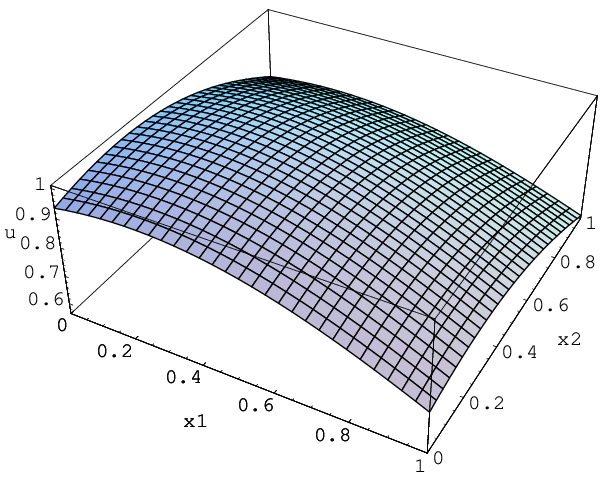}
\hfill
\includegraphics[width= 0.3 \columnwidth]{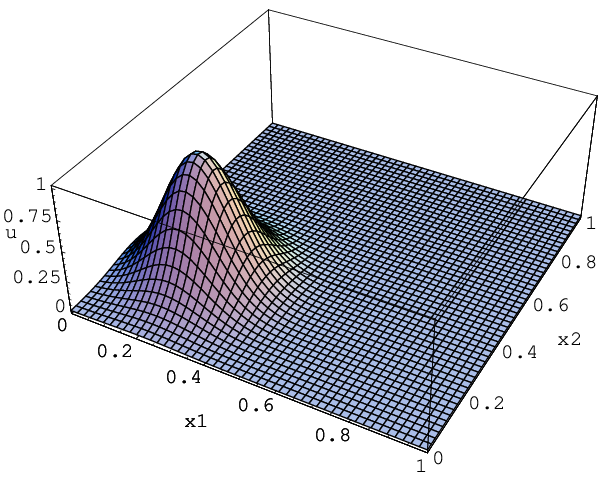}
\hfill
\includegraphics[width= 0.3 \columnwidth]{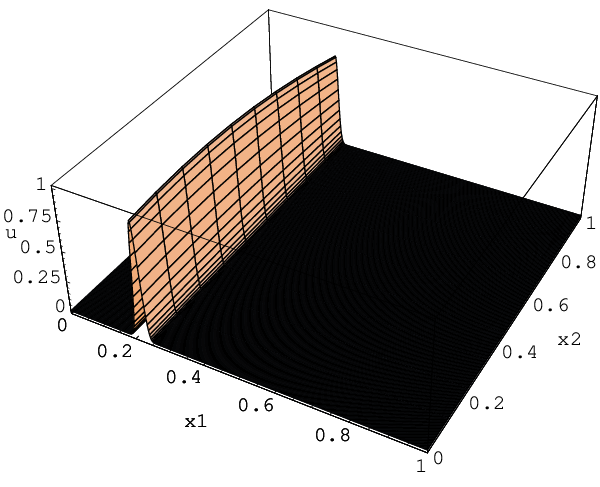}
\end{center}
\caption[Anisotropic solutions]{$u(x_1,x_2,p_3)$ for $(\lambda_1, \lambda_2, \lambda_3) = (1,1,1)$ (left),
$(\lambda_1, \lambda_2, \lambda_3) = (100,100,1)$ (middle) and
$(\lambda_1, \lambda_2, \lambda_3) = (10000,1,1)$ (right).}
\label{fig:lambdaPlots}
\end{figure}
\end{center}
Fig.~\ref{fig:lambdaPlotsConv} shows how the error increases with increasing (mixed)
derivatives.

\subsection{Parabolic and hyperbolic problems}

We consider
\[
u_t - \nu \Delta u + \vf b \cdot \nabla u = 0
\]
in two dimensions
for different values of $\nu$. The case $\nu = 0$ resembles the advection
equation.

\begin{figure}[ht]
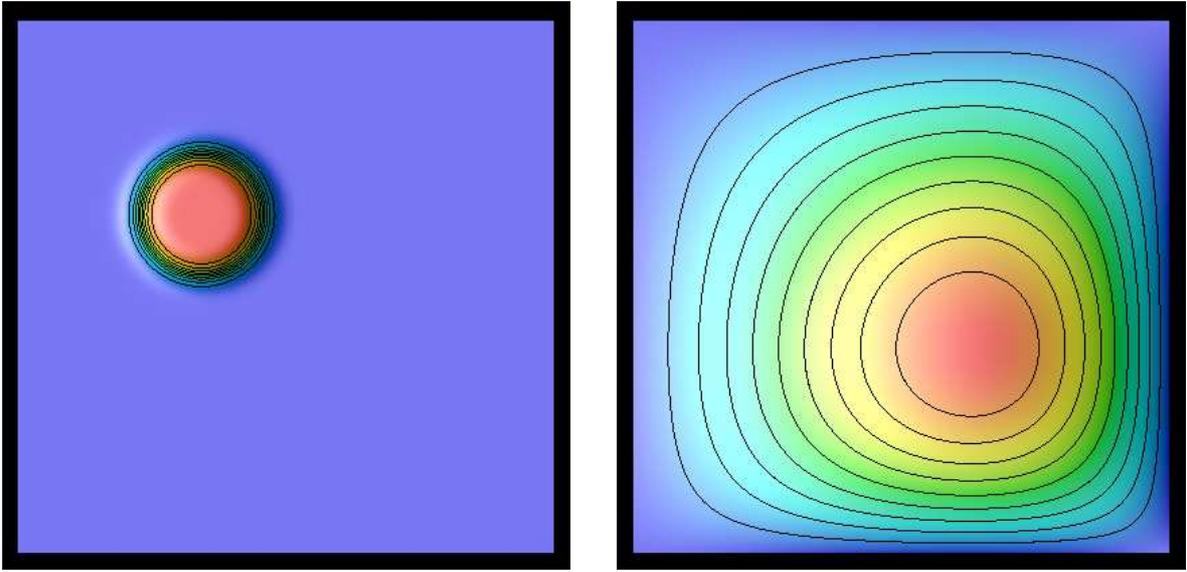

\begin{center}
\includegraphics[clip = true, trim = 65 170 65 170,angle =270,
width= 0.48 \columnwidth]{eps01nu09t0.epsf} \hfill
\includegraphics[clip = true, trim = 65 170 65 170,angle=270,
width= 0.48 \columnwidth]{eps01nu09.epsf}
\end{center}
\caption[Sparse grid solution for advection diffusion]{Sparse grid solution at level $n = 15$ 
for $\nu = 0.1, \epsilon = 0.9$ at $t=0$ (left) and $t=1$ (right).}
\label{fig:plot12}
\end{figure}

To adjust the smoothness, we choose an initial profile
\[
u(x_1,x_2,0) = \left\{\begin{array}{rl}
1 &\;\;\; r < \underline{r}, \\
\arctan \left(\tan^k \left( \frac{\pi}{2} \frac{|r-\underline{r}-\epsilon|}{\epsilon}\right)\right)
&\;\;\; \underline{r} \le r \le \overline{r}, \\
0 &\;\;\; r > \overline{r},
\end{array}\right.
\]
where
\begin{eqnarray*}
r = \sqrt{(x_1-m_1)^2 + (x_2-m_2)^2},
\end{eqnarray*}
$\underline{r} = (1-\epsilon) \overline{r}$ for $\epsilon \in [0,1]$
and $k\ge 1$ (see Fig.~\ref{fig:plot12}, left, for $\epsilon = 0.9$).

The value is 1 in an inner circle with centre $(m_1,m_2)$
and radius $\underline{r}$, and changes to 0 within a distance of $\epsilon$.
The regions are joined together such that $k-1$ is the order of
differentiability for $\epsilon>0$.
$\epsilon = 0$ is the discontinuous limit.
We choose $k=5$, because the analysis predicts that mixed derivatives
of order up to $2+2$ are required for upwinding on a sparse grid.
Homogeneous Dirichlet conditions $u=0$ are set at the boundary.
To avoid any
effects arising from the velocity being aligned with the grid,
we choose $b_1=0.31415926535897932385$, $b_2 = -0.27182818284590452354$,
we furthermore set $m_1 = 0.5 (1 - b_1)$, $m_2 = 0.5 (1 - b_2)$
and $\overline{r} = 0.5 |\vf b|_2$, such that in the non-diffusive case the profile
starts from the upper left quarter of the unit square, with the outer circle
going through the centre (Fig.~\ref{fig:plot12}, left, for $\epsilon = 0.9$),
and is propagated
down to the lower right quarter, touching the centre from the other side at $T=1$.
We evaluate the solution at the centre $x_* = (1/2,\ldots,1/2)$
and since the exact solution is unknown, we study the surplus
\begin{equation}
\label{pt-error2}
\hat{\epsilon}_n = |u_{n+1}(x_*)-u_n(x_*)|.
\end{equation}

\begin{center}
\begin{figure}
\psfrag{n}[l][r][1.0]{$n$}
\psfrag{delta}[l][u][1.0]{$\log_{2} \hat{\epsilon}_n$}
\begin{center}
\includegraphics[width= 0.6 \columnwidth]{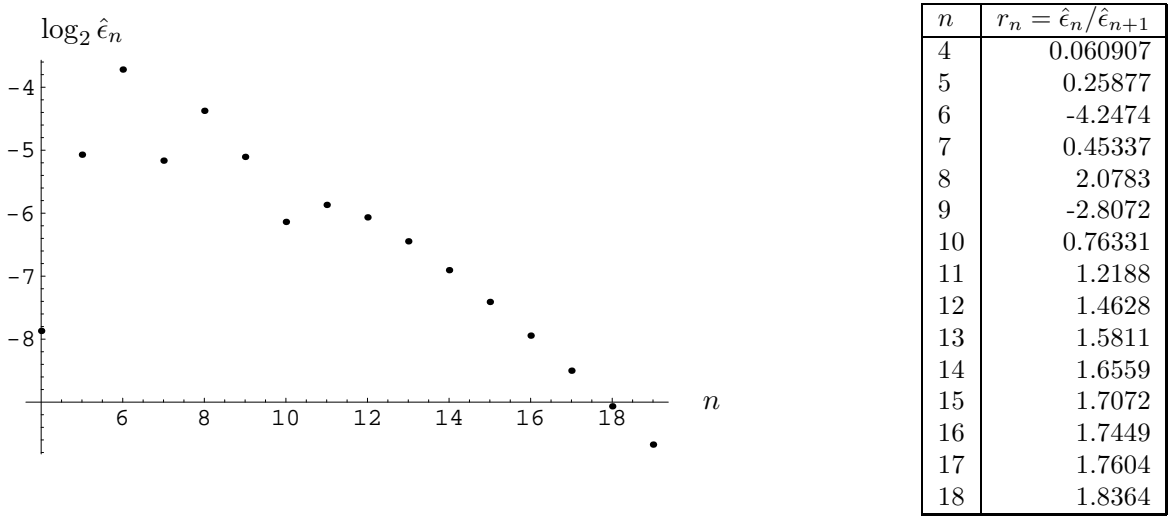}
\hfill
{\small
\begin{minipage}{0.22 \columnwidth}
\vspace{-5 cm}

\begin{tabular}{|l|r|}\hline
$n$ & $r_n = \hat{\epsilon}_n/\hat{\epsilon}_{n+1}$ \\ \hline
4 &     0.060907 \\
5 &     0.25877 \\
6 &     -4.2474 \\
7 &     0.45337 \\
8 &     2.0783 \\
9 &     -2.8072 \\
10 &    0.76331 \\
11 &    1.2188 \\
12 &    1.4628 \\
13 &    1.5811 \\
14 &    1.6559 \\
15 &    1.7072 \\
16 &    1.7449 \\
17 &    1.7604 \\
18 &    1.8364 \\ \hline
\end{tabular}
\vspace{0 cm}
\end{minipage}}
\end{center}
\caption[Convergence rate for advection diffusion]{$\hat{\epsilon}_n$ as in (\ref{pt-error2}) and convergence rate $r_n$
at refinement level $n$ for $\nu = 0.1$, $\epsilon = 0.9$.}
\label{fig:delta}
\end{figure}
\end{center}

We first consider the setting $\nu = 0.1$, $\epsilon = 0.9$. The numerical solution
is shown in Fig.~\ref{fig:plot12}.
First order convergence is observed as illustrated by Fig.~\ref{fig:delta}.

\begin{center}
\begin{figure}[ht]
\psfrag{n}[l][r][1.0]{$n$}
\psfrag{eps}[l][u][1.0]{$\log_{2} \hat{\epsilon}_n$}
\psfrag{plot1}[l][l][1.0]{$\nu = 0.1$}
\psfrag{plot2}[l][l][1.0]{$\nu = 0.01$}
\psfrag{plot3}[l][l][1.0]{$\nu = 0.001$}
\psfrag{plota}[l][l][1.0]{$\epsilon = 0.1$}
\psfrag{plotb}[l][l][1.0]{$\epsilon = 0.01$}
\psfrag{plotc}[l][l][1.0]{$\epsilon = 0.001$}
\psfrag{n0}[l][l][1.0]{$\nu = 0$}
\psfrag{e0}[l][l][1.0]{$\epsilon = 0$}
\begin{center}
\includegraphics[width= 0.48 \columnwidth]{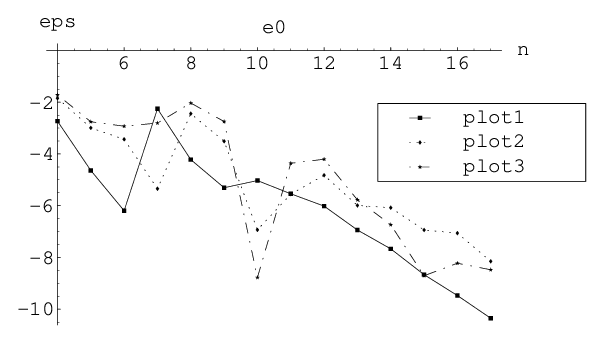} \hfill
\includegraphics[width= 0.48 \columnwidth]{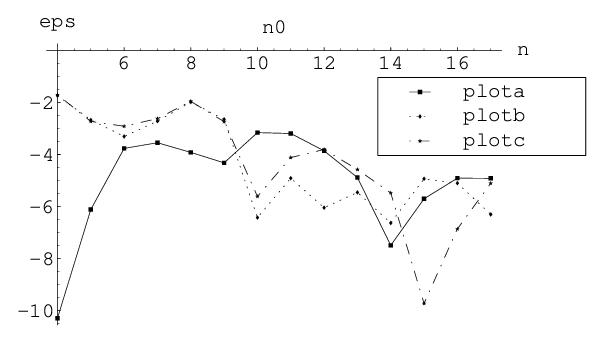}
\end{center}
\caption[Convergence for discontinuous data]{Convergence of the sparse grid solution for discontinuous
initial data for small diffusivity (left) and smooth initial data without
diffusion (right) for with $\hat{\epsilon}_n$ from (\ref{pt-error2}).}
\label{fig:deltax0}
\end{figure}
\end{center}

We now take $\epsilon = 0$ and change $\nu = 0.1, 0.01, 0.001$.
(Fig.~\ref{fig:deltax0}, left). Due to the discontinuous initial condition,
the problem does not fit into the theoretical framework. Nonetheless, the
smoothing property of the diffusion operator suffices to provide sufficient
regularity for $t>0$.

Alternatively,
taking $\nu = 0$, we let the smooth transition collapse to
$\epsilon = 0.1, 0.01, 0.001$ (Fig.~\ref{fig:deltax0}, right).
Although the problem is sufficiently smooth, the scales involved are not resolved
properly and the (mixed) derivatives are too large to reach the asymptotic range
at feasible refinement levels.

\begin{center}
\begin{figure}
\psfrag{n}[l][r][1.0]{$n$}
\psfrag{delta}[l][u][1.0]{$\log_{2} \hat{\epsilon}_n$}
\begin{center}
\includegraphics[width= 0.6 \columnwidth]{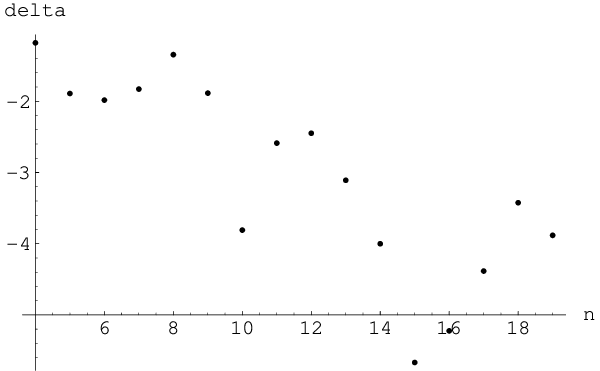}
\hfill
{\small
\begin{minipage}{0.22 \columnwidth}
\vspace{-5 cm}

\begin{tabular}{|l|r|}\hline
$n$ & $r_n = \hat{\epsilon}_n/\hat{\epsilon}_{n+1}$ \\ \hline
4 & 2.0381 \\
5 & 1.0974       \\
6 & 0.85633      \\
7 & -0.61664     \\
8 & 1.7166       \\
9 & -6.8473      \\
10 & 0.29471     \\
11 & 0.87042     \\
12 & 1.9361      \\
13 & 2.4390      \\
14 & -5.2997     \\
15 & -0.64137    \\
16 & -0.43160    \\
17 & 0.38301 \\
18 & 1.58094     \\ \hline
\end{tabular}
\vspace{0 cm}
\end{minipage}}
\end{center}
\caption[Convergence for discontinuous advection problem]{$\hat{\epsilon}_n$ as in (\ref{pt-error2}) and convergence rate $r_n$
at refinement level $n$ for discontinuous initial data without diffusion.}
\label{fig:delta4}
\end{figure}
\end{center}

Finally, consider the limiting case $\nu = 0$, $\epsilon = 0$, i.e.\ a discontinuous
initial profile without diffusion.
\begin{center}
\begin{figure}[ht]
\begin{center}
\begin{minipage}{0.48 \columnwidth}
\includegraphics[clip = true, trim = 65 170 65 170,angle=270,
width= \columnwidth]{eps0nu0t0.epsf}
\end{minipage}
\hfill
\begin{minipage}{0.48 \columnwidth}
\includegraphics[clip = true, trim = 65 170 65 170,angle=270,
width= \columnwidth]{eps0nu0.epsf}
\end{minipage}
\end{center}
\caption[Sparse grid solution for advection problem]{Sparse grid solution at level $n = 15$ for $\nu = 0, \epsilon = 0$
at $t=0$ (left) and $t=1$ (right).}
\label{fig:plot34}
\end{figure}
\end{center}
Already in Fig.~\ref{fig:plot34} the problems of the sparse grid become apparent.
It is noticeable how the anisotropic elements in the sparse grid fail to capture
the discontinuous transition at the circumference of the circle.
Fig.~\ref{fig:delta4} confirms that convergence (if any) is slow and erratic.


\section{Discussion}
\label{sec:discussion}

In this article, we give explicit error bounds for the sparse grid combination solution 
to model problems in a finite difference context,
and at the same time provide a framework that can be applied to more general
linear PDEs.
The ingredients are:
\begin{enumerate}
\item
sufficiently smooth and compatible data to yield solutions with bounded mixed derivatives
of required order;
\item
a discretisation scheme that provides a truncation error of certain mixed order;
\item
stability of the discretisation scheme, i.e.\ a bounded inverse in the maximum
norm.
\end{enumerate}

In most cases, the truncation error can be assessed easily by Taylor expansion.
The boundedness of the inverse will be harder to establish, but note that no
additional requirements to the corresponding full grid case are necessary in this regard.
Examples that fall into this category, e.g.\ where discrete maximum principles are
known, are elliptic and parabolic equations. The numerical results reproduce the
theoretical findings nicely.

Limitations we encountered concern non-smooth problems, dimensions in excess
of eight, and non-diagonal diffusion tensors. The latter can often be resolved
in practice by diagonalising the diffusion tensor and performing a principal
component or asymptotic analysis, see e.g.\ \citep{reisinger:06}.

{
\small
\begin{center}%
{\bfseries Acknowledgements}
\vspace{-.5em}%
\end{center}

The author wishes to thank Mike Giles for illuminating
discussions on the
subject, in particular for pointing him towards the concept of adjoint error
correction \citep{giles-pierce:01}, which proved vital for
the theory of Section \ref{sec:error-expansions};
Endre S{\"u}li and Tony Ware for advice on questions regarding
the regularity of solutions, which were raised by an anonymous referee of an earlier version of this manuscript whose
contribution in spotting this is also gratefully acknowledged.
}

\bibliographystyle{plainnat}

\bibliography{Literature}

\begin{thebibliography}{33}
\providecommand{\natexlab}[1]{#1}
\providecommand{\url}[1]{\texttt{#1}}
\expandafter\ifx\csname urlstyle\endcsname\relax
  \providecommand{\doi}[1]{doi: #1}\else
  \providecommand{\doi}{doi: \begingroup \urlstyle{rm}\Url}\fi

\bibitem[Bachmayr(2010)]{bachmayr:10}
M.~Bachmayr.
\newblock Hyperbolic wavelet discretization of the two-electron {S}chr\"odinger
  equation in an explicitly correlated formulation.
\newblock {AICES} Preprint 2010/06-2, RWTH Aachen, June 2010.

\bibitem[Bungartz(1992)]{bungartz:92}
H.-J. Bungartz.
\newblock \emph{D\"unne Gitter und deren Anwendung bei der adaptiven L\"osung
  der dreidimensionalen Poisson-Gleichung}.
\newblock PhD thesis, Technische Universit\"at M\"unchen, 1992.

\bibitem[Bungartz(1998)]{bungartz:98}
H.-J. Bungartz.
\newblock Finite elements of higher order on sparse grids.
\newblock Habilitationsschrift, Technische Universit\"at M\"unchen, 1998.

\bibitem[Bungartz and Griebel(2004)]{bungartz-griebel:04}
H.-J. Bungartz and M.~Griebel.
\newblock Sparse grids.
\newblock \emph{Acta Numerica}, 13:\penalty0 1--123, 2004.

\bibitem[Bungartz et~al.(1994)Bungartz, Griebel, R\"oschke, and
  Zenger]{bungartz-griebel-roeschke-zenger:94}
H.-J. Bungartz, M.~Griebel, D.~R\"oschke, and C.~Zenger.
\newblock Pointwise convergence of the combination technique for the {L}aplace
  equation.
\newblock \emph{East-West J. Num. Math.}, 2:\penalty0 21--45, 1994.

\bibitem[Dauge and Stevenson(2010)]{stevenson:10}
M.~Dauge and R.~Stevenson.
\newblock Sparse tensor product wavelet approximation of singular functions.
\newblock \emph{SIAM J. Math. Anal.}, 42\penalty0 (5):\penalty0 2203--2228,
  2010.

\bibitem[Dijkema et~al.(2009)Dijkema, Schwab, and Stevenson]{stevenson:09}
T.~J. Dijkema, C.~Schwab, and R.~Stevenson.
\newblock An adaptive wavelet method for solving high-dimensional elliptic
  {PDE}s.
\newblock \emph{Constr. Approx.}, 30\penalty0 (3):\penalty0 423--455, 2009.

\bibitem[Garcke and Griebel(2000)]{garcke-griebel:00}
J.~Garcke and M.~Griebel.
\newblock On the computation of the eigenproblems of hydrogen and helium in
  strong magnetic and electric fields with the sparse grid combination
  technique.
\newblock \emph{Journal of Computational Physics}, 165\penalty0 (2), 2000.

\bibitem[Gerstner and Griebel(2003)]{gerstner-griebel:03}
T.~Gerstner and M.~Griebel.
\newblock {D}imension--{A}daptive {T}ensor--{P}roduct {Q}uadrature.
\newblock \emph{Computing}, 71\penalty0 (1):\penalty0 65--87, 2003.

\bibitem[Gerstner et~al.(2009)Gerstner, Griebel, and
  Holtz]{gerstner-griebel-holtz:09}
T.~Gerstner, M.~Griebel, and M.~Holtz.
\newblock Efficient deterministic numerical simulation of stochastic
  asset-liability management models in life insurance.
\newblock \emph{Insurance: Math. Economics}, 44:\penalty0 434--446, 2009.

\bibitem[Giles and S\"uli(2002)]{giles-pierce:01}
M.B. Giles and E.~S\"uli.
\newblock Adjoint methods for {PDE}s: \emph{a posteriori} error analysis and
  postprocessing by duality.
\newblock \emph{Acta Numerica}, pages 145--236, 2002.

\bibitem[Griebel(2006)]{griebel06}
M.~Griebel.
\newblock Sparse grids and related approximation schemes for higher dimensional
  problems.
\newblock In L.-M. Pardo, A.~Pinkus, E.~S\"uli, and M.~Todd, editors,
  \emph{Foundations of Computational Mathematics 2005}, pages 106--161.
  Cambridge University Press, 2006.

\bibitem[Griebel and Hamaekers(2007)]{griebel-hamaekers:07}
M.~Griebel and J.~Hamaekers.
\newblock Sparse grids for the {S}chr\"odinger equation.
\newblock \emph{{ESAIM:} Mathematical Modelling and Numerical Analysis},
  41\penalty0 (2):\penalty0 215--247, 2007.

\bibitem[Griebel and Holtz(2010)]{griebelholtz}
M.~Griebel and M.~Holtz.
\newblock Dimension-wise integration of high-dimensional functions with
  applications to finance.
\newblock \emph{J. Complexity}, 26\penalty0 (5):\penalty0 455--489, 2010.

\bibitem[Griebel and Thurner(1995)]{griebel-thurner:95}
M.~Griebel and V.~Thurner.
\newblock The efficient solution of fluid dynamics problems by the combination
  technique.
\newblock \emph{Int. J. Num. Meth. for Heat and Fluid Flow}, 1995.

\bibitem[Griebel et~al.(1992)Griebel, Schneider, and
  Zenger]{griebel-schneider-zenger}
M.~Griebel, M.~Schneider, and C.~Zenger.
\newblock A combination technique for the solution of sparse grid problems.
\newblock In P.~de~Groen and R.~Beauwens, editors, \emph{Iterative Methods in
  Linear Algebra}. IMACS, Elsevier, North Holland, 1992.

\bibitem[Hegland(2003)]{Hegl2003}
M.~Hegland.
\newblock Adaptive sparse grids.
\newblock In K.~Burrage and Roger~B. Sidje, editors, \emph{Proc. of 10th
  Computational Techniques and Applications Conference CTAC-2001}, volume~44,
  pages 335--353, 2003.

\bibitem[Hilber et~al.(2005)Hilber, Matache, and
  Schwab]{hilber-matache-schwab:05}
N.~Hilber, A.M. Matache, and C.~Schwab.
\newblock Sparse wavelet methods for option pricing under stochastic
  volatility.
\newblock \emph{J. Comp. Fin.}, 8\penalty0 (4):\penalty0 1--42, 2005.

\bibitem[Kangro and Nicolaides(2001)]{kangro}
R.~Kangro and R.~Nicolaides.
\newblock Far field boundary conditions for {B}lack-{S}choles equations.
\newblock \emph{SIAM J. Numer. Anal.}, 38\penalty0 (4):\penalty0 1357--1368,
  2001.

\bibitem[Leentvaar and Oosterlee(2008)]{leentvaar-oosterlee:08}
C.C.W. Leentvaar and C.W. Oosterlee.
\newblock On coordinate transformation and grid stretching for sparse grid
  pricing of basket options.
\newblock \emph{J. Comp. Appl. Math.}, 222:\penalty0 193--209, 2008.

\bibitem[Nitsche(2005)]{nitsche}
P.-A. Nitsche.
\newblock Sparse approximation of singularity functions.
\newblock \emph{Constr. Approx.}, 21:\penalty0 63Ð81, 2005.

\bibitem[Novak and Ritter(1996)]{novak}
E.~Novak and K.~Ritter.
\newblock High dimensional integration of smooth functions over cubes.
\newblock \emph{Numer. Math.}, 75\penalty0 (1):\penalty0 79--98, 1996.

\bibitem[Pflaum(1997)]{pflaum:97}
C.~Pflaum.
\newblock Convergence of the combination technique for second-order elliptic
  differential equations.
\newblock \emph{SIAM J. Numer. Anal.}, 34\penalty0 (6):\penalty0 2431--2455,
  December 1997.

\bibitem[Pflaum and Zhou(1999)]{pflaum-zhou:99}
C.~Pflaum and A.~Zhou.
\newblock Error analysis of the combination technique.
\newblock \emph{Numerische Mathematik}, 84:\penalty0 327--350, December 1999.

\bibitem[Reisinger(2004)]{reisinger:04}
C.~Reisinger.
\newblock \emph{Efficient Numerical Methods for High-Dimensional Parabolic
  Equations and Applications in Option Pricing}.
\newblock PhD thesis, Universit\"at Heidelberg, 2004.

\bibitem[Reisinger and Wittum(2007)]{reisinger:06}
C.~Reisinger and G.~Wittum.
\newblock Efficient hierarchical approximation of high-dimensional option
  pricing problems.
\newblock \emph{SIAM J. Sci. Comp.}, 29\penalty0 (1), 2007.

\bibitem[Schwab et~al.(2008)Schwab, S\"uli, and Todor]{suli}
C.~Schwab, E.~S\"uli, and R.-A. Todor.
\newblock Sparse finite element approximation of high-dimensional
  transport-dominated diffusion problems.
\newblock \emph{{ESAIM:} Mathematical Modelling and Numerical Analysis},
  42\penalty0 (5):\penalty0 777--820, 2008.

\bibitem[Smolyak(1963)]{smolyak}
S.~Smolyak.
\newblock Quadrature and interpolation formulas for tensor products of certain
  classes of functions.
\newblock \emph{Soviet Math. Dokl.}, 4:\penalty0 240Ð243, 1963.

\bibitem[von Petersdorff and Schwab(2004)]{petersdorff-schwab:02}
T.~von Petersdorff and C.~Schwab.
\newblock Numerical solution of parabolic equations in high dimensions.
\newblock \emph{{ESAIM}: Mathematical Modelling and Numerical Analysis},
  38\penalty0 (1):\penalty0 93--127, 2004.

\bibitem[Wloka(1987)]{wloka:87}
J.~Wloka.
\newblock \emph{Partial differential equations}.
\newblock Cambridge University Press, 1987.

\bibitem[Yserentant(2010)]{yserentant:10}
H.~Yserentant.
\newblock \emph{Regularity and Approximability of Electronic Wave Functions}.
\newblock Number 2000 in Lecture Notes in Mathematics. Springer, 2010.

\bibitem[Zeiser(2010)]{zeiser:10}
A.~Zeiser.
\newblock Wavelet approximation in weighted {S}obolev spaces of mixed order
  with applications to the electronic {S}chr\"odinger equation.
\newblock {SPP} 1324 {P}reprint, TU Berlin, July 2010.

\bibitem[Zenger(1990)]{zenger:90}
C.~Zenger.
\newblock Sparse grids.
\newblock In W.~Hackbusch, editor, \emph{Parallel Algorithms for Partial
  Differential Equations}, volume~31 of \emph{Notes on Numerical Fluid
  Dynamics}, 1990.
\newblock Proceedings of the Sixth GAMM-Seminar.

\end{thebibliography}

\begin{appendix}

\section{Smoothness of solutions to the Poisson problem on the hyper-cube}\label{app:smoothness}

\newcommand{\dtwo}[1]{\partial_{#1}^2}
\newcommand{\dfour}[1]{\partial_{#1}^4}
\newcommand{\dk}[1]{\partial_{#1}^k}

\subsection{Fourier series expansion}
\label{app:fourier}

In this section, we derive a sufficient condition on the decay of the
Fourier (sine) coefficients of the right-hand side $f$, in order for
the solution $u$ to the Poisson problem (\ref{poisson}) to have sufficient
regularity.

\begin{theorem}

Let $f_{\vf k}$ be the coefficients of $f$ in a sine-series with multi-index
$\vf k = (k_1,\ldots, k_d) \in \mathbb{N}^d$. Then if
\[
\sum_{\vf k} \sum_{\vf i, |\vf i|_1\le \lceil d/2\rceil}
\frac{k_1^{8+i_1} \cdot \ldots \cdot k_d^{8+i_d}}{\left(k_1^2+\ldots+k_d^2\right)^2}
f_{\vf k}^2 < \infty,
\]
the solution $u$ to the Poisson problem has regularity
\[
u \in X_4^d
\]
with $X_4^d$ defined in (\ref{eqn:space})
as the space of all functions with continuous mixed derivatives
up to order 4 that vanish at the boundaries.
\end{theorem}

\begin{proof}
Expanding both $f$ and $u$ into $\sin$-series,
\begin{eqnarray*}
f(\vf x) &=& \sum_{\vf k} f_{\vf k} \prod_{j=1}^d \sin(\pi k_j x_j) \\
u(\vf x) &=& \sum_{\vf k} u_{\vf k} \prod_{j=1}^d \sin(\pi k_j x_j)
\end{eqnarray*}
then for the Fourier coefficients $u_{\vf k}$ and $f_{\vf k}$
\[
u_{\vf k} = \frac{f_{\vf k}}{|\vf k|_2^2}.
\]
Norm equivalence gives
\begin{eqnarray*}
f \in L_2 &\Leftrightarrow& \sum_{\vf k} f_{\vf k}^2 < \infty \\
\frac{\partial^{4 d} f}{\partial x_1^4 \ldots \partial x_d^4}
\in L_2 &\Leftrightarrow&
\sum_{\vf k} \left[k_1^4 \cdot \ldots \cdot k_d^4 f_{\vf k}\right]^2 < \infty \\
\frac{\partial^{4 d} u}{\partial x_1^4 \ldots \partial x_d^4}
\in L_2 &\Leftrightarrow&
\sum_{\vf k} \left[\frac{k_1^4 \cdot \ldots \cdot k_d^4}{k_1^2+\ldots+k_d^2}
\right]^2 f_{\vf k}^2 < \infty
\end{eqnarray*}

Further, see e.g.\ \citep{wloka:87}, Theorem 6.2, p.~107,
$W_2^k \hookrightarrow C^l$ for $k-l\ge d/2$, in particular
\[
W_2^k \hookrightarrow C \text{ for } k\ge d/2.
\]
From this the statement follows.
\end{proof}

\subsection{Maximum norm of mixed derivatives}
\label{app:derivatives}

\begin{theorem}
For the Poisson problem (\ref{poisson}),
if continuous derivatives up to the relevant order exist, then
\[
\|\dfour{k_1} \ldots \dfour{k_m} u\|_\infty 
\le \frac{1}{8} \|\dfour{k_1} \ldots \dfour{k_m} f\|_\infty 
+ \sum_{i=1}^d \|\dfour{k_1} \ldots \dfour{k_m}\dtwo{i} f\|_\infty.
\]

\end{theorem}

\begin{proof}
In order to get bounds on the derivative norms, we derive PDEs and boundary
conditions for mixed (fourth) derivatives, by differentiating the PDE,
and can then apply a maximum principle argument.
From $\sum_{i=1}^d \dtwo{i} u = f$
follows
\begin{equation}
\label{eqn:mixed}
\sum_{i=1}^d \dtwo{i} \left( \dfour{k_1}\ldots \dfour{k_d} u\right) =
\dfour{k_1}\ldots \dfour{k_d} f,
\end{equation}
and therefore derivatives of the solution satisfy a Poisson problem with
the right-hand side a derivative of the original one. We now derive boundary conditions
for $\dfour{k_1}\ldots \dfour{k_d} u$
at $x_1 = 0$. The other cases follow by permutation and symmetry.
From
\[
\dk{i} u = 0 \text{ for } i>1,\; k\ge 0
\]
follows
\[
\dtwo{1} u = \dtwo{1} u + \sum_{i=2}^d \dtwo{i} u = f,
\]
and from $\dtwo{i} \dtwo{j} u =0$ for $i,j>1$
\[
\dtwo{1} \dtwo{j} u = \dtwo{1} \dtwo{j} u + \sum_{i=2}^d \dtwo{i} \dtwo{j} u
= \dtwo{j} f
\]
for $j>1$.
Also,
\[
\dfour{1} u + \sum_{i=2}^d \dtwo{1} \dtwo{i} u = \dtwo{1} f,
\]
and consequently for $2\le k_1,\ldots,k_m\le d$
\begin{equation}
\label{eqn:bc-mixed}
\dfour{k_1} \ldots \dfour{k_m} \dfour{1} u = \dfour{k_1} \ldots \dfour{k_m} \dtwo{1} f
- \sum_{i=2}^d \dfour{k_1} \ldots \dfour{k_m} \dtwo{i} f
\end{equation}
for $x_1=0$. All other mixed derivatives without $x_1$ in them are zero at
this boundary face.

From (\ref{eqn:mixed}) and (\ref{eqn:bc-mixed}) the result follows by a
maximum principle argument.
\end{proof}

\section{Proofs from Section \ref{sec:interpolation}}
\label{app:proofs3}

\begin{proof}[Proof of Lemma \ref{lemma-expansion-high}]
Define
\begin{eqnarray*}
v &:=&
u(x_1,\ldots,x_d) - u(0,\ldots,0) - \sum_{i=1}^d x_i \, \partial_i u(x-x^{(i)})
\; - \sum_{\scriptsize\begin{array}{c}i,j=1\\ i \neq
    j\end{array}}^d x_i x_j \, \partial_i \partial_j u(x-x^{(i,j)})
- \ldots - \\
&& \hspace{-0.9 cm}
- \; x_1 \cdot \ldots \cdot x_d \left.\partial_1 \ldots \partial_d u
\right|_{x_1, \ldots, x_d = 0}  \; - \;
\sum_{i=1}^d \int_0^{x_i} (x_i-s_i) \, \partial_i^2 u(s^{(i)}) \ds_i - \ldots
- \\
&& \hspace{-0.9 cm} - \hspace{-0.5 cm}
\sum_{\scriptsize \begin{array}{c}\{i_1,\ldots,i_{d-1}\}\\\subset
    \{1,\ldots,d\}\end{array}} \!\!
\int_0^{x_{i_1}} \ldots \int_0^{x_{i_{d-1}}} \! (x_{i_1}-s_{i_1})\cdot\ldots\cdot
(x_{i_{d-1}}-s_{i_{d-1}})
\partial_{i_1}^2 \ldots \partial_{i_{d-1}}^2 u(s^{(i_1,\ldots,i_{d-1})})
\, \ds_{i_{d-1}} 
\ldots \ds_{i_{1}}.
\end{eqnarray*}
%
It is straightforward to see that 
\[
v(x_1,\ldots,x_d) = 0 \quad {\rm if} \quad \exists k: \; x_k = 0
\]
and also
\[
D^{\boldsymbol \alpha} v(x_1,\ldots,x_d) = 0
\quad {\rm if} \quad \exists k: \; x_k = 0 \; \wedge \; \alpha_k = 1.
\]
From
\[
\partial_1^2 \ldots \partial_d^2 v = 0
\]
one gets $v=0$ from which the result follows.
\end{proof}

\begin{proof}[Proof of Theorem \ref{lemma-approx}]
Without loss of generality consider a point $x$ in the box
$B:=[0,h_1]\times \ldots\times [0,h_d]$.
For all corner points $p = (p_1,\ldots,p_d) \in
B_c := \bigotimes_{i=1}^d \{0,h_i\}$
use Lemma \ref{lemma-expansion-high} to express
%
\begin{eqnarray}
\nonumber
u(p) &=& u(x) + \sum_{i=1}^d (p_i-x_i)\partial_i u(x^{(i)})
+ \sum_{\scriptsize\begin{array}{c}i,j=1\\ i \neq
    j\end{array}}^d (p_i-x_i) (p_j-x_j) \partial_i \partial_j u(x^{(i,j)})
+ \ldots + \vspace{-1 cm} \\
\nonumber
&& \hspace{-1 cm} + \; (p_1-x_1) \cdot \ldots \cdot (p_d-x_d) 
\partial_1 \ldots \partial_d u(x) \; + \;
\sum_{i=1}^d \int_{x_i}^{p_i} (p_i-s_i) \partial_i^2 u(s^{(i)}) \ds_i +
\ldots + \\
&& \hspace{-1 cm} + \; \int_{x_1}^{p_1} \ldots \int_{x_d}^{p_d}  
(p_1-s_1)\cdot\ldots\cdot (p_d-s_d) \,
\partial_1^2 \ldots \partial_d^2 u(s) \ds_d \ldots \ds_1
\label{exp-origin}
\end{eqnarray}
Here
$x=(x_1,\ldots,x_d)$, and introduce
$x^{(i_1,\ldots,i_k)}, s^{(i_1,\ldots,i_k)}\in \R^d$
with $x^{(i_1,\ldots,i_k)}_j = x_j$ if $j\in \{i_1,\ldots i_k\}$, $p_j$
otherwise, and $s^{(i_1,\ldots,i_k)}_j = s_j$ if $j\in \{i_1,\ldots i_k\}$,
$x_j$ otherwise.

Then insert in the multilinear approximation
\begin{eqnarray*}
(\mathcal I u(\vf x_h))(x) = \sum_{p \in B_c} 
\frac{|p_1-x_1|}{h_1} \cdot \ldots \cdot
\frac{|p_d-x_d|}{h_d} \cdot u(p)
= u(x) + \sum_{m=1}^d
\!\! \sum_ {
\tiny \begin{array}{c}
\{i_1,\ldots,i_m\} \\
\subset \{1,\ldots,d\}
\end{array}}
\!\! \sum_{p \in B_c} w_{i_1,\ldots,i_m}(p,x)
\end{eqnarray*}
with
%
\begin{eqnarray*}
w_{i_1,\ldots,i_m}(p,x)= \left(\prod_{k=1}^d \frac{|p_k-x_k|}{h_k}\right)
\int_{x_{i_1}}^{p_{i_1}} \!\!\ldots
\int_{x_{i_m}}^{p_{i_m}}
\!\left( \prod_{k=1}^m (p_{i_k}-s_{i_k}) \right)
\, \partial_{i_1}^2 \ldots \partial_{i_m}^2 u(s^{(i_1,\ldots,i_m)}) \ds_{i_m}
\ldots \ds_{i_1},
\end{eqnarray*}
because all the terms in (\ref{exp-origin}) that are multilinear are exactly
represented.
By inserting the points $p \in B_c$ one
shows after some calculation the representation 
\begin{equation}
\label{interpol-alpha}
\alpha_{j_1,\ldots,j_m} (x;h_{j_1},\ldots, h_{j_m}) =
\prod_{k=1}^m \frac{x_{j_k}}{h_{j_k}} \left(1- \frac{x_{j_k}}{h_{j_k}} \right)
\!\! \sum_{\scriptsize \begin{array}{c} \vf y, \vf z \in \{0,1\}^m \\ y_k + z_k = 1 \\
k=1,\ldots,m \end{array}}
\prod_{k=1}^m
\left(\frac{x_{j_k}}{h_{j_k}}\right)^{\! y_k} \!\left(1- \frac{x_{j_k}}{h_{j_k}}
\right)^{\! z_k} w_{\vf y,\vf z}(x;\vf h)
\end{equation}
with $\vf h=(h_{i_1},\ldots,h_{i_m})$ and
%
\[
w_{\vf y, \vf z}(x; \vf h) = \frac{1}{\prod\limits_{k=1}^d x_k^{2 y_k}
  (h_{j_k}-x_k)^{2 z_k}} \int_{x_1}^{z_1 h_{j_1}} \!\! \ldots
\int_{x_m}^{z_m h_{j_m}}
\left(\prod_{k=1}^m (z_k h_{j_k}-s_{k})\right)
\partial_1^2 \ldots \partial_m^2 u(s^{(j_1,\ldots,j_m)})
\ds_m \ldots \ds_1.
\]
The rest follows because the maximum of $\xi^2 (1-\xi)$ on $[0,1]$ is
$\frac{4}{27}$.
\end{proof}

\begin{proof}[Proof of Theorem \ref{theorem:error-expansion-total}]
From Corollary \ref{corollary-smoothness}
$w_{i_1,\ldots,i_k}(\vf x_h; h_{i_1},\ldots,h_{i_k})$ (in (\ref{rep-expansion}))
is the restriction of a function $w_{i_1,\ldots,i_k}$ 
from hyperplanes $\vf I_h^{(i_1,\ldots,i_k)}$ to the grid $\vf x_h$, where
for the derivatives in the continuous directions
\[
\|D^{\boldsymbol \alpha} w_{i_1,\ldots,i_k}(\cdot; h_{i_1},\ldots,h_{i_k})\|_0
\le \frac{k!}{48^k} \|D^{\boldsymbol \beta} u \|_\infty
\]
with $\alpha_i = 4$, $i \notin \{i_1,\ldots,i_k\}$ $\beta_i = 4$ holds.

By Theorem \ref{lemma-approx}
for points $x$ on the hyper-planes $\vf I_h^{(i_1,\ldots,i_k)}$
\begin{eqnarray}
\nonumber
(\mathcal I w_{i_1,\ldots,i_k}(\vf x_h; h_{i_1},\ldots,h_{i_k}))(x) &=& 
w_{i_1,\ldots,i_k}(x; h_{i_1},\ldots,h_{i_k}) + \\
&& \hspace{-4 cm}
\sum_{\scriptsize \begin{array}{c}\{j_1,\ldots, j_m\} \\
\cap \{i_1,\ldots,i_k\} = \emptyset \end{array}} 
\gamma_{i_1,\ldots, i_k;j_1,\ldots,j_m}
(x,h_{i_1},\ldots, h_{i_k}; h_{j_1},\ldots, h_{j_m}) h_{j_1}^2\cdot \ldots \cdot
h_{j_m}^2
\label{expansion-interpol}
\end{eqnarray}
with
$
|\gamma_{i_1,\ldots, i_k;j_1,\ldots,j_m}
(x;h_{i_1},\ldots, h_{i_k}; h_{j_1},\ldots, h_{j_m})|
\le \frac{4^m}{27^m} \|D^{\boldsymbol \alpha} w\|_\infty
$
where $\alpha_{i_l} = 4, \alpha_{j_n} = 2$.
Values between the hyperplanes are obtained by multilinear interpolation, which
introduces no new maxima, and thus
\begin{eqnarray}
\nonumber
\mathcal I (u(\vf x_h) - \vf u_h) &=& \!\!
\sum_{\scriptsize \begin{array}{c}\{i_1,\ldots, i_k\} \\
\subset \{1,\ldots,k\} \end{array}} \!\!
\beta_{i_1,\ldots,i_k}(\cdot; h_{i_1},\ldots,h_{i_k})
h_{i_1}^2\cdot \ldots \cdot h_{i_k}^2 + \\
&& \hspace{-3 cm} + \!\! \sum_{\scriptsize \begin{array}{c}\{j_1,\ldots, j_m\} \\
\cap \{i_1,\ldots,i_k\} = \emptyset \end{array}}\!\!\!\!\!\!\!\!
\beta_{i_1,\ldots,i_k;j_1,\ldots,j_m}(\cdot; h_{i_1},\ldots,h_{i_k};
 h_{j_1},\ldots,h_{j_m})
h_{i_1}^2\cdot \ldots \cdot h_{i_k}^2 \cdot
h_{j_1}^2\cdot \ldots \cdot h_{j_m}^2
\label{interp}
\end{eqnarray}
with
$
\|\beta_{i_1,\ldots, i_k;j_1,\ldots,j_m}
(\cdot;h_{i_1},\ldots, h_{i_k}; h_{j_1},\ldots, h_{j_m})\|_\infty
\le K \frac{k!}{96^k} \frac{4^m}{27^m}.
$

Putting (\ref{interp}) and (\ref{error-expansion-approx}) together one gets
\[
u - \mathcal{I} \vf u_h = \sum_{m=1}^{d}
\sum_{
\tiny \begin{array}{c}
\{j_1,\ldots,j_m\} \\
\subset \{1,\ldots,d\}
\end{array}}
v_{j_1,\ldots,j_m}
(\cdot;h_{j_1},\ldots, h_{j_m}) h_{j_1}^2\cdot \ldots \cdot h_{j_m}^2
\]
with
\begin{eqnarray*}
\|v_{j_1,\ldots,j_m}
(\cdot; h_{j_1},\ldots, h_{j_m})\|_\infty
&\le&
K \frac{4^m}{27^m} + 
K \frac{m!}{96^m} +
K \sum_{l=1}^{m-1} \bino{m}{l}\left(\frac{4}{27}\right)^{l}
\frac{(m-l)!}{96^{(m-l)}} \\
&=& K \sum_{l=0}^{m} \bino{m}{l}\left(\frac{4}{27}\right)^{l}
 \frac{(m-l)!}{96^{m-l}}
= K \frac{m!}{96^m}\sum_{l=0}^{m} \frac{1}{l!} \left(\frac{384}{27}\right)^l \\
&<& K \frac{m!}{96^m} \e^{384/27} < 150188 \, K \frac{m!}{96^m}.
\end{eqnarray*}

\end{proof}

\begin{proof}[Proof of Theorem \ref{theorem:error-expansion-total-upw}]
Follows by combining Theorems \ref{theorem:expansion-on-grid-adv} and
\ref{lemma-approx} similarly to the Poisson
case by observing that one can write second order terms (from the bilinear
interpolation) as first order terms by defining
\[
\beta_{i_1,\ldots,i_k}(\cdot;h_{i_1},\ldots,h_{i_k}) h_{i_1}^2 \cdot
\ldots \cdot
h_{i_k}^2 =: \hat{\beta}_{i_1,\ldots,i_k}(\cdot;h_{i_1},\ldots,h_{i_k}) h_{i_1}
\cdot \ldots \cdot h_{i_k}
\]
with $|\hat{\beta}_{i_1,\ldots,i_k}|\le |\beta_{i_1,\ldots,i_k}|$ and the
final estimate becomes
\[
\|v_{j_1,\ldots,j_m}
(\cdot; h_{j_1},\ldots, h_{j_m})\|_\infty \le 
K \frac{m!}{2^m} \sum_{l=0}^m \frac{1}{l!} \left(\frac{8}{27}\right)^l \le K
\frac{m!}{2^m} \e^{8/27} \le \frac{3}{2} K \frac{m!}{2^m}.
\]
\end{proof}

\section{Proofs from Section \ref{sec:error-bounds}}
\label{app:proofs4}

\begin{proof}[Proof of Lemma \ref{err-rep}]
Inserting gives
\begin{eqnarray*}
F(n) = \lefteqn{
\sum_{|\vf{i}|=n} v(2^{-i_1},\ldots,2^{-i_m})
2^{-p i_1} \cdot \ldots \cdot 2^{-p i_m} =} \hspace{2 cm} \\
&=& 
\sum_{|\vf{i}|=n} v(2^{-i_1},\ldots,2^{-i_m})
2^{-p \sum_{k=1}^m i_k} \\
&=& \sum_{l=0}^n 
\underbrace{
\sum_{\sum_{k=1}^m i_k = l}
v(2^{-i_1},\ldots,2^{-i_m})}_{=:s_l}
2^{-pl}
\sum_{\sum_{k=m+1}^d i_k = n-l} 1 \\
&=& \sum_{l=0}^n s_l 2^{-pl} \bino{n-l+d-m-1}{d-m-1}.
\end{eqnarray*}
From Lemma \ref{tech} below, \ref{alltech2}.\ ($d\rightarrow d-m$,
$f_l\rightarrow s_l 2^{-pl}$) one obtains
\begin{equation}
\label{dieFormel}
\dlt^{d-m} F(n) =
\dlt^{d-m} \sum_{|\vf{i}|=n}
v(2^{-i_1},\ldots,2^{-i_m})
2^{-p i_1}\cdot \ldots \cdot 2^{-p i_m} = s_{n+d-m} 2^{-p(n+d-m)}.
\end{equation}
For all $j\ge 0$
\[
\dlt^j 2^{-pn} = \left(2^{-p}-1\right)^j 2^{-pn},
\]
as one sees inductively ($j\rightarrow j+1$) by
\[
\dlt^{j+1} 2^{-pn} = \left(2^{-p}-1\right)^j
\left(2^{-p(n+1)}-2^{-pn}\right) = \left(2^{-p}-1\right)^{j+1} 2^{-pn}.
\]
Therefore it follows from Proposition \ref{prod-diff} that
\begin{eqnarray*}
\lefteqn{
\dlt^{m-1} \left[s_{n+d-m} 2^{-p(n+d-m)}\right] =} \hspace{2 cm} \\
&=& \sum_{j=0}^{m-1} \bino{m-1}{j} \dlt^j s_{n+d-m} \dlt^{m-1-j}
2^{-p(n+d-m+j)}\\
&=& \sum_{j=0}^{m-1} \bino{m-1}{j} \dlt^j s_{n+d-m}
\left(2^{-p}-1\right)^{m-1-j} 2^{-p(n+d-m+j)}\\
&=& 2^{-p(n+d-1)} \sum_{j=0}^{m-1} \bino{m-1}{j} \left(1-2^p\right)^{m-1-j}
\dlt^j s_{n+d-m} \\
&=& 2^{-p(n+d-1)} \left(\dlt-2^p+1\right)^{m-1} s_{n+d-m} \\
&=& 2^{-p(n+d-1)} \sum_{i=0}^{m-1} s_{n+d-m+i} \bino{m-1}{i}
(-2)^{p (m-1-i)} \\
&=& 2^{-p (n+d-1)} \sum_{i=0}^{m-1} s_{n+d-i-1} \bino{m-1}{i}
(-2)^{p i}.
\end{eqnarray*}
\end{proof}

 \begin{lemma}[Differencing formula]
\label{tech}
Let $d \in \N$ and $f, F \in \R^{\N_0}$ s.~t.
\[
F(n) = \sum_{l=0}^n f_l \bino{n-l+d-1}{d-1}.
\]
\begin{enumerate}
\item
\label{alltech1}
Then for $0\le k< d$
\[
\dlt^k F(n) = G^k(n) + H^k(n),
\]
where
\begin{eqnarray*}
G^k(n) &:=& \left\{
\begin{array}{ll}
0 & k = 0 \\
\sum_{j=1}^k f_{n+j} \bino{d-j-1}{k-j} & k\ge 1 
\end{array}
\right., \\
H^k(n) &:=&
\sum_{l=0}^n f_{l} \bino{n-l+d-1}{d-k-1}.
\end{eqnarray*}
\item
\label{alltech2}
\[
\dlt^{d} F(n) = f_{n+d}
\]
\end{enumerate}
\end{lemma}
\begin{proof}
\begin{enumerate}
\item
Induction in $k$:
Clearly ($k=0$)
\[
\dlt^0 F(n) = F(n) = H^0(n) = G^0(n) + H^0(n).
\]
Since $F(n)$ has the form $F(n) = \sum_{l=0}^n f(l,n)$,
\[
\dlt F(n) = f(n+1,n+1) -
\sum_{l=0}^n \left(f(l,n+1)-f(l,n)\right),
\]
hence
\begin{eqnarray*}
\dlt F(n) &=& f_{n+1} + \sum_{l=0}^n
f_l \left[\bino{n-l+d}{d-1}-\bino{n-l+d-1}{d-1} \right] \\
&=& f_{n+1} + \sum_{l=0}^n f_l \bino{n-l+d-1}{d} \\
&=& G^1(n) + H^1(n),
\end{eqnarray*}
i.~e.\ the case $k=1$.

For $k>1$ the contributions ($k \rightarrow k+1$) are
\begin{eqnarray*}
\dlt G^k(n) &=& \sum_{j=1}^k 
\left(f_{n+j+1}-f_{n+j} \right)
\bino{d-j-1}{k-j} \\
&\hspace{-2 cm}=& \hspace{-1 cm} f_{n+k+1} - f_{n+1} {\small \bino{d-2}{k-1}}
+ \sum_{j=2}^k f_{n+j} \left[ {\small
\bino{d-j}{k-j+1}-\bino{d-j-1}{k-j} }\right] \\
&\hspace{-2 cm}=& \hspace{-1 cm}\sum_{j=2}^{k+1} f_{n+j}
\bino{d-j-1}{k-j+1} - f_{n+1}\bino{d-2}{k-1} \\
&\hspace{-2 cm}=& \hspace{-1 cm}G^{k+1}(n) - f_{n+1}\bino{d-2}{k} 
- f_{n+1} \bino{d-2}{k-1} \\
&\hspace{-2 cm}=&\hspace{-1 cm} G^{k+1}(n) - f_{n+1} \bino{d-1}{k}
\end{eqnarray*}
and in $H^k(n)$
\begin{eqnarray*}
\dlt H^k(n) &=& f_{n+1} {\small\bino{d-1}{d-k-1}} + \sum_{l=0}^n
f_l {\small \left[\bino{n-l+d}{d-k-1}-\bino{n-l+d-1}{d-k-1}\right]}
\\
&=& f_{n+1} \bino{d-1}{k} + \sum_{l=0}^n
f_l \bino{n-l+d-1}{d-k-2} \\
&=& H^{k+1}(n) + f_{n+1} \bino{d-1}{k}.
\end{eqnarray*}
This completes the result 
as
\begin{eqnarray*}
\dlt^{k+1} F(n) &=& \dlt \dlt^k F(n) = \dlt \left(G^k(n) + H^k(n)\right) =
\dlt G^k(n) + \dlt H^k(n) \\
&=& G^{k+1}(n) + H^{k+1}(n).
\end{eqnarray*}
\item
From \ref{alltech1}.\ one obtains for $k=d-1$
\[
\dlt^{d-1} F(n) = \sum_{j=1}^{d-1} f_{n+j} + \sum_{l=0}^n f_l
\]
and
\[
\dlt^d F(n) = \sum_{j=1}^{d-1} (f_{n+j+1}-f_{n+j}) +
f_{n+1} = f_{n+d}.
\]
\end{enumerate}
\end{proof}

\begin{proof}[Proof of Theorem \ref{zweiteOrdnung}]
Because of Lemma \ref{konsistenz} (the exact solution $u$ is constant over
the grid levels)
\[
\dlt^{d-1} \sum_{|\vf i| = n} u = u \; \dlt^{d-1} \! \sum_{|\vf i| = n} 1
= u \; \dlt^{d-1} N(n,d) = u,
\]
and one may write
\[
u-u_n = \dlt^{d-1} \sum_{|\vf{i}|=n} (u-U(\vf i)).
\]
The error terms in
\[
u-U(\vf i) = \sum_{m=1}^{d} \sum_{\scriptsize \begin{array}{c}\{j_1,\ldots, j_m\} \\
\subset \{1,\ldots,d\} \end{array}} v_{j_1,\ldots,j_m}
(\cdot;2^{-i_{j_1}},\ldots, 2^{-i_{j_m}}) 2^{-p i_{j_1}}
\cdot \ldots \cdot 2^{-p i_{j_m}}
\]
will now be studied separately with the help of Lemma \ref{err-rep}, applied to
\[
F_{j_1,\ldots,j_m}(n) := \sum_{|\vf{i}|=n} v_{j_1,\ldots,j_m}
(\cdot,2^{-i_{j_1}},\ldots,2^{-i_{j_m}})
2^{-p i_{j_1}} \cdot \ldots \cdot 2^{-p i_{j_m}}
\]
in
\begin{equation}
\label{compact-formula}
u-u_n = \sum_{m=1}^d \sum_{\scriptsize
  \begin{array}{c}\{j_1,\ldots, j_m\} \\ \subset \{1,\ldots,d\} \end{array}}
\dlt^{d-1} F_{j_1,\ldots,j_m} (n).
\end{equation}
The point is that with Lemma \ref{err-rep} the factors $h_i$ then no longer
appear separately, but only in their highest order $h_1 \cdot \ldots \cdot h_n
= 2^{-n}$.
It remains to estimate the coefficients, which leads to the polynomial terms
in $n$.
From $|v_{j_1,\ldots,j_m}|\le K$ follows
\[
|s_l| \le K \bino{l+m-1}{m-1} \quad \Rightarrow \quad
\max_{i=0}^{m-1}|s_{n+d-i-1}| \le K \bino{n+d+m-2}{m-1},
\]
therefore
\begin{eqnarray}
\label{absEst}
\bigg|\sum_{i=0}^{m-1} s_{n+d-i-1} \bino{m-1}{i}
(-2)^{pi}\bigg|
&\le& 
K \bino{n+d+m-2}{m-1} \sum_{i=0}^{m-1} \bino{m-1}{i} 2^{pi} \\
&\le& K \bino{n+d+m-2}{m-1} (2^p+1)^{m-1},
\nonumber
\end{eqnarray}
and with Lemma \ref{err-rep}
\begin{eqnarray}
\hspace{-1 cm}
\big|\dlt^{d-1} F_{j_1,\ldots,j_m}(n)
 \big|
\label{ersteUngl}
&\le&
2^{-p(n+d-1)} (2^p+1)^{m-1} K \bino{n+d+m-2}{m-1} \\
\label{zweiteUngl}
&\le&
2^{-p(n+d-1)} (2^p+1)^{d-1} K \bino{n+2(d-1)}{d-1} \\
\label{dritteUngl}
&<&
2^{-pn} \left(\frac{2^p+1}{2^p}\right)^{d-1} \!\!\! \frac{K}{(d-1)!} (n+2(d-1))^{d-1}.
\end{eqnarray}
Since the number of terms in (\ref{compact-formula}) is $2^d-1$, we get
\[
|u-u_n| < \left(\frac{2^p+1}{2^{p-1}}\right)^{d-1} \!\!\! \frac{2 K}{(d-1)!} (n+2(d-1))^{d-1} 2^{-pn}.
\]
\end{proof}

\begin{proof}[Proof of Lemma \ref{supiAbsch}]
A binomial of the form (\ref{ersteUngl}) or
(\ref{zweiteUngl}), respectively, can be written as ($l=n+d-1$, $k=m-1$)
\begin{eqnarray*}
\bino{l+k}{k} = \frac{\prod_{j=1}^k (l+j)}{\prod_{j=1}^k j}
= \prod_{j=1}^k (1+\frac{l}{j}).
\end{eqnarray*}
From the inequality between the arithmetic and geometric mean one gets
\[
\sqrt[k]{\prod_{j=1}^k (1+\frac{l}{j})} < \frac{1}{k}\left[k + l \left(\sum_{j=1}^k
  \frac{1}{j} \right)\right],
\]
furthermore
\begin{eqnarray*}
\sum_{j=1}^k \frac{1}{j} &<& 1 + \int_1^k \frac{\dx}{x} = 1 + \ln k,
\end{eqnarray*}
because the first sum is a lower sum for the integral.
This gives
\[
\bino{l+k}{k} < \left[1+l \frac{1 + \ln k}{k}\right]^k
\]
and
\[
\bino{n+d+m-2}{m-1} < \left[1+(n+d-1) \frac{1 + \ln (m-1)}{m-1}\right]^{m-1}.
\]
The rest follows as in the proof of Theorem \ref{zweiteOrdnung}.
\end{proof}

\begin{proof}[Proof of Lemma \ref{exactAsympt}]
For continuous $v_{1,\dots,m}$ we can first show
\begin{equation}
\label{limit}
\lim_{n\rightarrow \infty} \frac{s_n}{\bino{n+m-1}{m-1}} =
v_{1,\ldots,m}(0,\ldots,0) =: v_0.
\end{equation}
In the following the index of $v_{1,\ldots,m}$ is omitted for simplicity
of notation.
Let $\epsilon>0$ and
\[
c := \sup_{0\le h_1,\ldots,h_m \le 1} |v(h_1,\ldots,h_m)-v_0|.
\]
Choose $n_0$ such that
\[
\forall i_k \ge n_0, 1\le k \le m: \quad
|v(2^{-i_1},\ldots,2^{-i_m})-v_0| \le \frac{\epsilon}{2}
\]
and then $n$ sufficiently large such that
\[
N_0 := \bino{n-k n_0+m-1}{m-1} \ge \left(1-\frac{\epsilon}{c}\right)
\underbrace{\bino{n+m-1}{m-1}}_{=:N}.
\]
The latter is possible, because
\[
\bino{n-k n_0+m-1}{m-1} \Big/ \bino{n+m-1}{m-1} = \prod_{j=1}^{m-1}
\left(1-\frac{k n_0}{n+j} \right) \rightarrow 1
\]
for $n\rightarrow \infty$.
The idea is now to show that the contribution of the terms that do not lie in
an $\epsilon$-ball around $v_0$ can be neglected.
This motivates the splitting
\[
s_n = \sum_{\scriptsize \begin{array}{c}\sum_{k} i_k =n \\ \forall k: i_k \ge n_0
\end{array}}
v(2^{-i_1},\ldots,2^{-i_m})
 \;\;\; +
\sum_{\scriptsize \begin{array}{c}\sum_{k} i_k =n \\ \exists k: i_k < n_0
\end{array}}
v(2^{-i_1},\ldots,2^{-i_m}).
\]
Because of
\[
\sum_{\scriptsize \begin{array}{c}\sum_{k} i_k =n \\ \forall k: i_k \ge n_0
\end{array}} 1
=
\sum_{\scriptsize \begin{array}{c}\sum_{k} i_k =n-m n_0 \\ \forall k: i_k \ge 0
\end{array}} 1
= N_0
\]

one sees
\begin{eqnarray*}
s_n - v_0 N =
\bigg[
\sum_{\scriptsize
  \begin{array}{c}\sum_{k} i_k =n \\ \forall k: i_k \ge n_0
\end{array}} 
+
\sum_{\scriptsize
  \begin{array}{c}\sum_{k} i_k =n \\ \exists k: i_k < n_0
\end{array}}\bigg]
\!\!
\left(v(2^{-i_1},\ldots,2^{-i_m})-v_0\right)
\end{eqnarray*}

and consequently
\begin{eqnarray*}
\left|s_n - v_0 N\right| \le 
\frac{\epsilon}{2} N + c \left(N-N_0\right)
\le \frac{\epsilon}{2} N + c \frac{\epsilon}{2c} N
= \epsilon N.
\end{eqnarray*}
Division by $N$ leads to (\ref{limit}).

Asymptotically one gets instead of (\ref{absEst})
for $v_0 \neq 0$
\begin{eqnarray*}
\sum_{i=0}^{m-1} s_{n+d-i-1} \bino{m-1}{i}
(-2)^{pi}
&\sim& v_0 N \sum_{i=0}^{m-1} \bino{m-1}{i} (-2)^{p(m-1-i)} \\
&=&
v_0 \bino{n+m-1}{m-1} (1-2^p)^{m-1},
\end{eqnarray*}
in other words $\forall \epsilon>0 \;\, \exists N \;\, \forall n\ge N$
\begin{eqnarray}
\nonumber
\bigg|\sum_{i=0}^{m-1} s_{n+d-i-1} \bino{m-1}{i}
(-2)^{pi}\bigg| \le v_0 \bino{n+m-1}{m-1}
(2^p-1)^{m-1} (1+\epsilon).
\label{epsinequ}
\end{eqnarray}
The rest follows as in Theorem \ref{zweiteOrdnung}.
\end{proof}

\end{appendix}

\newpage
\listoftables
\listoffigures

\end{document}